\documentclass{article}
\usepackage[a4paper,top=4cm,bottom=4cm,left=2.2cm,right=2.2cm,bindingoffset=5mm]{geometry}

\usepackage{xcolor}
\newcommand{\Rone}[1]{\textcolor{black}{#1}}
\newcommand{\Rtwo}[1]{\textcolor{black}{#1}}

\usepackage[hyphens]{url}
\usepackage[breaklinks=true,colorlinks,linkcolor={blue},citecolor={red},urlcolor={purple},]{hyperref}

\usepackage{todonotes}

\usepackage{cite}

\newcommand\Item[1][]{%
  \ifx\relax#1\relax  \item \else \item[#1] \fi
  \abovedisplayskip=0pt\abovedisplayshortskip=0pt~\vspace*{-\baselineskip}}

\usepackage{subcaption}
\usepackage{graphicx}
\usepackage{caption}
\usepackage{enumitem}
\usepackage{multirow}
\usepackage{pict2e,overpic}
\usepackage{pgfplots}

\graphicspath{{figure/}{dati/burgers/}{dati/burgers/}{dati/convtest/}{dati/lintra/}{dati/systems/}{dati/wave_analysis/}}

\usepackage{amssymb,bm,amsmath,amsthm,amsfonts,accents}
\usepackage{graphicx,mathtools}	
\usetikzlibrary{calc}
\usepackage{esvect}
\usepackage{algorithm}
\usepackage{algorithmic}

\usepackage{pgfplotstable}
\usepackage{booktabs}

\newcommand{\INPUT}{\item[\textbf{input:}]}

\mathtoolsset{showonlyrefs}
			 
\newcommand{\disp}{\displaystyle}     
\newcommand{\R}{\mathbb R}
\newcommand{\C}{\mathbb C}
\newcommand{\Id}{\mathrm{I}}
\newcommand\norm[1]{\left\lVert#1\right\rVert}        
        
\def\dx{h}

\def\dt{\Delta t}

\newcommand{\kk}{\mathrm{k}}
\newcommand{\Fnum}{\mathcal{F}}

\newcommand{\DIRK}[1][]{\ensuremath{\mathsf{DIRK#1}}}
\newcommand{\DG}{{\mathsf{DG}}}
\newcommand{\RK}{{\mathsf{RK}}}

\newcommand{\DGRKdue}{\DG_2\RK_{2}}
\newcommand{\DGRKtre}{\DG_3\RK_{3}}
\newcommand{\DIRKdue}{\DG_2\DIRK_{22}^{0.25}}
\newcommand{\DIRKlstab}{\DG_2\DIRK_{22}^\mathsf{Lstab}}
\newcommand{\DIRKlstabTre}{\DG_{3}\DIRK_{33}^\mathsf{Lstab}}
\newcommand{\DIRKlam}[1]{\DG_2\DIRK_{22}^{#1}}

\newcommand{\DIRKQuattroTre}{\DG_3\DIRK_{43}}

\newcommand{\Imag}{\mathfrak{Im}}
\newcommand{\Real}{\mathfrak{Re}}

\let\oldparagraph=\paragraph
\renewcommand\paragraph[1]{\oldparagraph{#1.}}

\theoremstyle{plain}
\newtheorem{defn}{Definition}[section]

\newtheorem{prop}[defn]{Proposition}

\newtheorem{example}[defn]{Example}
\newtheorem{remark}[defn]{Remark}
\numberwithin{equation}{section}

\numberwithin{defn}{section}
\usepackage{upref}

\title{\Large\textbf{Dissipation-dispersion analysis of fully-discrete implicit discontinuous Galerkin methods and application to stiff hyperbolic problems}}

\author{\normalsize{Maya Briani}\thanks{Istituto per le Applicazioni del Calcolo ``Mauro Picone'', Consiglio Nazionale delle Ricerche, Via dei Taurini 19, 00185 Rome, Italy (\href{mailto:maya.briani@cnr.it}{maya.briani@cnr.it}).}
\and {\setcounter{footnote}{3}\normalsize{Gabriella Puppo}\thanks{Dipartimento di Matematica ``Guido Castelnuovo'', Sapienza Universit\`a di Roma, P.le Aldo Moro 5, 00185 Rome, Italy (\href{mailto:gabriella.puppo@uniroma1.it}{gabriella.puppo@uniroma1.it}).}}
\and {\setcounter{footnote}{2}\normalsize{Giuseppe Visconti}\thanks{Dipartimento di Matematica ``Guido Castelnuovo'', Sapienza Universit\`a di Roma, P.le Aldo Moro 5, 00185 Rome, Italy (\href{mailto:giuseppe.visconti@uniroma1.it}{giuseppe.visconti@uniroma1.it}).}}
}
\date{\today}

\begin{document}
\allowdisplaybreaks
\maketitle

\begin{abstract}
	The application of discontinuous Galerkin (DG) schemes to hyperbolic systems of conservation laws requires a careful interplay between space discretization, carried out with local polynomials and numerical fluxes at inter-cells, and time-integration to yield the final update. An important concern is how the scheme modifies the solution through the notions of numerical dissipation-dispersion. As far as we know, no analysis of these artifacts has been considered for implicit integration of DG methods. The first part of this work intends to fill this gap, showing that the choice of the implicit Runge-Kutta impacts deeply on the quality of the solution. We analyze one-dimensional dissipation-dispersion to select the best combination of the space-time discretization for high Courant numbers.
	
	Then, we apply our findings to the integration of one-dimensional stiff hyperbolic systems. Implicit schemes leverage superior stability properties enabling the selection of time-steps based solely on accuracy requirements. High-order schemes require the introduction of local space limiters which make the whole implicit scheme highly nonlinear. To mitigate the numerical complexity, we propose to use appropriate space limiters that can be precomputed on a first-order prediction of the solution. Numerical experiments explore the performance of this technique on scalar equations and systems.
\end{abstract}

\paragraph{AMS subject classification} 65M60, 65L04, 35L65
\paragraph{Keywords} implicit methods, discontinuous Galerkin schemes, hyperbolic systems, diffusion-dispersion analysis.

\tableofcontents

\section{Introduction}

Hyperbolic partial differential equations (PDEs) are fundamental in the mathematical modeling of various physical phenomena. These equations describe the propagation of waves and signals at finite speeds in fluid dynamics, electromagnetism, geophysics, acoustics, to mention just a few applications.

The accurate and efficient numerical solution of hyperbolic problems is essential for predicting and analyzing these physical phenomena. One of the major challenges in solving hyperbolic PDEs is addressing stiffness, particularly in problems with multiple scales or stiff source terms~\cite{Mengaldo2019}. In fact, stiffness typically arises when there are rapid changes in the solution over very short time scales or when there is a large disparity in the magnitudes of the eigenvalues of the system's Jacobian matrix. There are many practical applications where the phenomenon of interest involves low-speed movement. For example, low-Mach number problems occur when the fluid velocity is much smaller than the speed of sound, leading to stiffness in the governing equations. In these cases, the fluid can be approximated as nearly incompressible, making the accurate numerical propagation of sound waves less critical. The literature on this topic is huge. We refer to~\cite{Abbate2017,Boscarino2018,Degond2011,Dellacherie2010,Dimarco2017,Tavelli2017,ThomannIolloPuppo}, and references therein.

From a numerical viewpoint, stiffness imposes severe restrictions on time-step sizes when using explicit time-integration methods, which significantly increases computational costs. In fact, the Courant-Friedrichs-Lewy (CFL) condition forces the time-step to be small enough to ensure stability, which in turn imposes accuracy on fast-moving waves. To better address the dynamics of the slower-moving fluid, implicit time-integration methods are often used. These methods are unconditionally stable and do not adhere to the CFL condition, allowing to choose the time-step based on the fluid's velocity rather than on the speed of sound. Indeed, this fact enables to concentrate computational power on the waves of interest, see~\cite{Puppo2024}, increasing their accuracy. \Rone{It is also important to note that, while implicit schemes enable larger time steps and accurate resolution of slow dynamics, this may come at the expense of reduced accuracy for fast acoustic waves when they are present.}

\begin{figure}[t!]
	\centering
	\includegraphics[width=\textwidth]{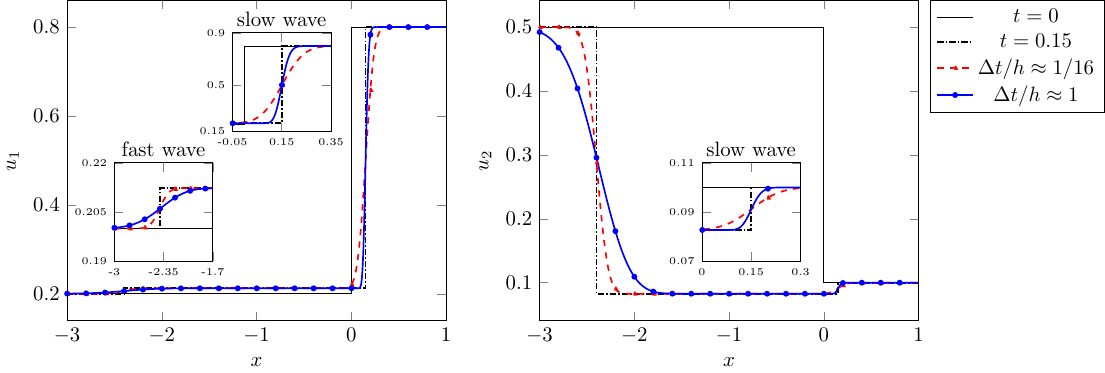}
	\caption{Numerical solution of a $2\times 2$ linear hyperbolic system obtained with the first-order backward Euler scheme. The dashed red line represents the approximation with a time-step computed on the fast wave, governed by the largest eigenvalue $\lambda_1 \approx -16$. The solid blue line shows the approximation computed with a time-step based on the slow wave, governed by the smallest eigenvalue $\lambda_2 \approx 1$.\label{fig:linear:system}}
\end{figure}

\Rtwo{We refer to Figure~\ref{fig:linear:system} to illustrate, on a simple model problem, the interplay between time-step selection and the accuracy on different wave families. In this toy example, we consider a $2\times 2$ linear hyperbolic system with eigenvalues $\lambda_1 \approx -16$ and $\lambda_2 \approx 1$. Using a first-order backward Euler scheme, we compare two choices of the time-step: $\Delta t_{\text{fast}} = h/16$, based on the fast wave (dashed red line), and $\Delta t_{\text{slow}} = h$, based on the slow wave (solid blue line). Here, $h$ represents the mesh-size. As shown, the smaller time-step $\Delta t_{\text{fast}}$ (which an explicit scheme would be constrained to use) provides resolution of the fast wave but is overly diffusive on the slow one. Conversely, the larger step $\Delta t_{\text{slow}}$, feasible only with an implicit scheme, is much more accurate on the slow wave.}

\Rtwo{This simple illustration highlights a more general principle: implicit time integration enables the user to shift the balance of accuracy towards selected wave families, at the price of solving nonlinear algebraic systems at each step. Hence, while the example is deliberately simplified, the underlying trade-off between stability, accuracy, and computational cost is representative of broader classes of hyperbolic problems.}

One class of implicit methods that has gained attention for hyperbolic problems is the Diagonally Implicit Runge-Kutta (DIRK) method~\cite{CarpenterKennedy,Ketcheson2009}. DIRK methods offer a good balance between stability and computational efficiency, making them a viable option for large-scale simulations of stiff hyperbolic systems.

In this work, we focus on the combination of DIRK time-integration methods with Discontinuous Galerkin (DG) space approximations.

The initial DG method incorporating Runge-Kutta time discretizations (RKDG) for solving hyperbolic PDEs was introduced by Cockburn and Shu~\cite{CS1989}. This approach was subsequently generalized and widely applied, see the review~\cite{Cockburn2005} and references therein. Error analysis for RKDG methods has been conducted by several researchers, including Zhang and Shu~\cite{Zhang2004, Zhang2006, Zhang2010}, Cockburn and Guzmán~\cite{Cockburn2008}, and Zhong and Shu~\cite{Zhong2011}.

DG methods employ local polynomial spaces within each cell. The use of high-degree polynomial approximations allows to achieve high-order accuracy, which is crucial for resolving fine-scale features in the solution. Unlike traditional finite element methods, the solution in DG methods is allowed to be discontinuous across cell boundaries. This feature is particularly advantageous for capturing sharp gradients and discontinuities, such as shock waves and contact discontinuities, which are common in hyperbolic problems. We refer to~\cite{Cockburn2001,Cockburn2003} for reviews. However, when using high-order methods, spurious oscillations can develop around discontinuities.
DG utilizes carefully designed space limiters within each element to control oscillations and maintain the non-oscillatory nature of the solution near discontinuities. Several approaches have been proposed~\cite{Krivodonova2004,Qiu2005,Krivodonova2007,BiswasDevineFlaherty94,Rider2001,Dumbser2014,Kuzmin2010,Kuzmin2014,Ortleb2020}, but a general consensus has not yet been reached, see for instance the review~\cite{Zhu2013}.

In this paper, first of all, we analyze one-dimensional dissipation and dispersion properties of DIRK-DG schemes. Dissipation-dispersion relations are essential tools for understanding how accurately the method captures wave propagation. Wave analysis helps in designing schemes that minimize numerical artifacts such as spurious oscillations and excessive damping. The dissipation-dispersion properties of numerical methods have been extensively studied, including DG methods. However, much of this research has focused on semi-discrete DG schemes~\cite{Ainsworth2004,Ainsworth2006,Hu2002,Hu1999} or fully-discrete explicit RKDG methods~\cite{Yang2013,Guo2013}.

To our knowledge, dissipation and dispersion analysis for fully-discrete implicit DG has not yet been studied.
To fill this gap, we derive analytical formulations of the dissipation and dispersion errors for second-order DIRK-DG schemes in terms of the CFL number for $\kk h \ll 1$, where $\kk$ is the wave number and $h$ the mesh-size.
The analysis of diffusion and dispersion for fully-discrete DG schemes increases tremendously in complexity with the order of the scheme. For this reason, we carry out this task analytically for second-order schemes, while we show a numerical approximation for third-order schemes.
The analysis demonstrates that the DIRK-DG scheme generates a principal, or physical, mode that closely approximates the exact dissipation-dispersion relations over a range of wave numbers, while other modes introduce spurious effects. This finding aligns with results from semi-discrete DG methods. This analysis provides valuable insights into optimizing the coupling between temporal and spatial discretizations to achieve desirable dissipation and dispersion characteristics.

We use insights of the proposed Fourier analysis, typically applied to linear equations, in order to design DIRK-DG schemes for nonlinear problems as well. However, implicit methods introduce nonlinearity through space limiters, which can be challenging to handle cost-effectively. To address this, we adapt a concept from~\cite{Puppo2023,Puppo2024} to design a framework that mitigates the nonlinearity introduced by space limiters, retaining only the unavoidable nonlinearity of the physical flux function. Our approach utilizes a first-order implicit predictor for pre-computing space limiters. Although first-order implicit schemes are fully linear with respect to the data on linear problems, they tend to produce significant dissipation errors. Thus, a high-order corrector is needed to improve resolution.

The implicit method proposed in~\cite{Puppo2023,Puppo2024} was specifically designed for a third-order finite-volume approximation of hyperbolic systems. This method achieved third-order accuracy by employing a third-order DIRK scheme for time integration and a third-order Central Weighted Essentially Non-Oscillatory (CWENO) reconstruction for spatial discretization, as detailed in~\cite{Cravero2018,Semplice2023}. The first-order implicit scheme used was based on a composite backward Euler method, evaluated at the DIRK abscissae, combined with a piecewise constant (linear in data) spatial reconstruction. This predictor allows to freeze the nonlinear weights of the CWENO reconstruction, resulting in a third-order implicit scheme where nonlinearity arises solely from the flux function. In this paper, we adopt the same predictor and adapt the scheme within the framework of DG methods.

We mention the following works on the implicit integration of hyperbolic systems. A similar idea of~\cite{Puppo2023,Puppo2024} was introduced in~\cite{Gottlieb2006} where an explicit predictor was used within a fifth-order implicit WENO scheme. Fully nonlinear implicit scheme have been more investigated. For instance, a third order RADAU time integrator combined with a third-order WENO reconstruction can be found in~\cite{Arbogast2020}. Fully implicit approaches for DG methods have been considered in~\cite{Southworth2022a,Southworth2022b}. \Rone{We also mention the work~\cite{Qin2018}, where a backward Euler DG scheme combined with a positivity-preserving strategy is shown to be effective in controlling oscillations near shocks. This robustness, however, is mainly due to the strong dissipative nature of backward Euler, which distinguishes it from the low-dissipation high-order DIRK schemes considered in the present study.}
Semi-implicit, implicit-explicit, local time-stepping, semi-Lagrangian and active flux treatments of stiff hyperbolic equations were also investigated, e.g., in~\cite{CNPT:10a,FZ:23,FKRZ:22,ABIR:19,BB:24,BoscarinoQiuRussoXiong,Despres23}.

The rest of the paper is organized as follows. In Section~\ref{sec:DG:time:approximation}, we provide a comprehensive review of the DG space discretization methods for hyperbolic systems, with a particular emphasis on implicit time-integration techniques. Section~\ref{sec:analysis} presents a detailed Fourier analysis of fully-discrete implicit DG schemes, where we define measures of dissipation and dispersion, focusing specifically on DIRK-DG schemes of second- and third-order. This section also includes analytical expressions for the eigenvalues and relative energies associated with second-order schemes. In Section~\ref{sec:limiting}, we introduce the concept of slope limiters in DG schemes and present the novel framework aimed at reducing nonlinearity by precomputing limiters using first-order implicit predictors. Finally, Section~\ref{sec:numerical:simulations} showcases one-dimensional numerical simulations for both scalar equations and systems, demonstrating the practical application and effectiveness of the proposed methods, which reflects the insights of the dissipation-dispersion analysis. We conclude the paper in Section~\ref{sec:conclusion} by summarizing the results and discussing future perspectives.

\section{Diagonally implict Runge-Kutta discontinuous Galerkin approximation} \label{sec:DG:time:approximation}

In this section we present the general formulation of Discontinuous Galerkin (DG) methods coupled with Diagonally Implicit Runge-Kutta (DIRK) time approximations applied to systems of $m>1$ conservation laws
\begin{equation} \label{eq:1dconsLaw}
\left\{\begin{array}{l}
	\mathbf{u}_t + \mathbf{f}_x(\mathbf{u}) = \mathbf{0} 
	\medskip\\
	\mathbf{u}(x,0)=\mathbf{u}_0(x)
\end{array}\right.
\end{equation}
where, $\mathbf{u}: (x,t) \in \mathbb{R} \times \mathbb{R}^+_0 \mapsto \mathbf{u}(x,t) \in \mathbb{R}^m$ is the quantity of interest, and $\mathbf{f}: \mathbb{R}^m \to \mathbb{R}^m$ is the vector of the flux functions.

\subsection{Space approximation}

First, we divide the computational domain $[0,L]$ into $N$ cells with boundary points
\begin{equation*}
0 = x_{1/2} < x_{3/2}<\cdots<x_{N+1/2}=L,
\end{equation*}
and we denote the center of cell $I_j=[x_{j-1/2},x_{j+1/2}]$ by $x_j$ and the length of cell $I_j$ by $\dx_j$. For the sake of simplicity, throughout the work we consider a uniform grid, thus $\dx_j=h, \ \forall j$. Then, we consider a polynomial approximation of the solution on each cell. To this end, we assume that the test function space is given by $V_h^p = \{v : v_{| I_j} \in \mathbb{P}_p (I_j)\}$, where $\mathbb{P}_p(I_j)$ is the space of polynomials of degree at most $p$ on cell $I_j$. As the local basis functions, we adopt the Legendre polynomials
\begin{equation}\label{eq:leg_pol}
\zeta_0(x) = 1, \quad \zeta_l(x) = \disp\frac{1}{2^l\, l!}\disp\frac{\mathrm{d}^l(x^2-1)^l}{\mathrm{d}x^l}, \quad l=1,\cdots,p.
\end{equation}
However, the procedure described below does not depend on the specific basis chosen for the polynomials. 
The approximate solution $\mathbf{u}_h$ is then expressed as follows:
\begin{equation}\label{eq:uh}
\mathbf{u}_h(x,t) = \sum_{l=0}^p \mathbf{u}^{(l)}_j(t)\, v^{(l)}_j(x) \quad \mbox{for } x\in I_j,
\end{equation}
where
\begin{equation}\label{eq:v_base}
v^{(l)}_j(x) = \zeta_l(2(x-x_j)/h)
\end{equation}
and $\mathbf{u}^{(l)}_j(t) \in \mathbb{R}^m$, $l = 0,\cdots, p$ are the \textit{moments} of the vector solution $\mathbf{u}(x,t)$ within the element $I_j$ and represent the degrees of freedom of the polynomial approximation.

To determine the approximate solution, multiply equation \eqref{eq:1dconsLaw} by test functions $v^{(l)}_j(x)$, $l = 0,\cdots, p$, then integrate over cell $I_j$, perform integration by parts, and thus, evolve the degrees of freedom $\mathbf{u}^{(l)}_j$ by the following system of equations for $l = 0,\ldots, p$:
\begin{equation}\label{eq:DG_odesystem}
		\disp\frac{h}{2l+1}\disp\frac{\mathrm{d}\mathbf{u}^{(l)}_j(t)}{\mathrm{d}t} = Q\left(\mathbf{u}_h(\cdot,t);l\right)
		- \left[ \Fnum(\mathbf{u}^-_{j+1/2}(t),\mathbf{u}^+_{j+1/2}(t)) - (-1)^l \Fnum(\mathbf{u}^-_{j-1/2}(t),\mathbf{u}^+_{j-1/2}(t)) \right],
\end{equation}
where $\mathbf{u}^\pm_{j+1/2}(t) = \mathbf{u}_h(x^\pm_{j+1/2},t)$ are the left and right limits of the discontinuous solution $\mathbf{u}_h$ at the cell interface $x_{j+1/2}$ and
\begin{equation} \label{eq:integral:discr}
	Q\left(\mathbf{u}_h(\cdot,t);l\right) \approx \disp\int_{-1}^1 \mathbf{f}(\mathbf{u}_h(x_j + \frac{h}{2}y,t))\disp\frac{\mathrm{d} \zeta_l(y)}{\mathrm{d}y} \mathrm{d}y, 
\end{equation}
for instance $Q$ can be a Gauss-Legendre quadrature rule of order matching the order of the DG approximation.
Finally, $\Fnum$ is a consistent and monotone (non-decreasing in the first argument and non-increasing in the second argument) numerical flux function that approximates the physical flux function at the interfaces, namely
\begin{subequations} \label{eq:numfluxfnc}
	\begin{equation}
		(\mathbf{u}^-,\mathbf{u}^+) \in \R^m\times\R^m \mapsto \Fnum(\mathbf{u}^-,\mathbf{u}^+) \in \R^m,
	\end{equation}
	such that
	\begin{equation}
		\Fnum\left(\mathbf{u}^-_{j+1/2}(t),\mathbf{u}^+_{j+1/2}(t)\right) \approx \mathbf{f}\left(\mathbf{u}_h(x_{j+1/2},t)\right).
	\end{equation}
\end{subequations}
Here, the function $\Fnum$ is any approximate Riemann solver, which is applied component-wise on vector-valued inputs. Numerical fluxes introduce the problem of the boundary conditions since $\mathcal{F}$ has to be evaluated on both sides of the domain boundaries. However, unlike finite-volume schemes, where boundary conditions often require information from ghost cells of an extended domain beyond the physical boundaries, DG methods can handle boundary conditions in a more localized manner.

\begin{example}[$\mathbb{P}_{p\leq2}$ DG approximation]
	In this work we will consider third-order schemes at most, namely approximations in space with $p\leq2$. The example below considers $p=2$. For lower order, $p<2$, you just drop the higher order modes and reduce the accuracy of the quadrature.\\
	The Legendre basis functions are
		\begin{align}
		\{ \zeta_0(x) = 1 \}, & \quad \text{for } \ p=0 \\
		\{\zeta_0(x), \zeta_1(x) = x\}, & \quad \text{for } \ p=1 \\
		\{\zeta_0(x), \zeta_1(x), \zeta_2(x) = \frac{1}{2} (3x^2-1)\}, & \quad \text{for } \ p=2
		\end{align}
	and the reconstructions of the solution at the cell interfaces become
	\begin{align*}
		\mathbf{u}_{j+1/2}^-(t) &= \sum_{l=0}^p \mathbf{u}_j^{(l)}(t) v_j^{(l)}(x_{j+1/2}) = \sum_{l=0}^p \mathbf{u}_j^{(l)}(t) \zeta_l(1) = \sum_{l=0}^p \mathbf{u}_j^{(l)}(t), \\
		\mathbf{u}_{j+1/2}^+(t) &= \sum_{l=0}^p \mathbf{u}_{j+1}^{(l)}(t) v_{j+1}^{(l)}(x_{j+1/2}) = \sum_{l=0}^p \mathbf{u}_{j+1}^{(l)}(t) \zeta_l(-1) = \sum_{l=0}^p (-1)^l \mathbf{u}_{j+1}^{(l)}(t).
 	\end{align*}
 	We use the three-point Gauss-Legendre quadrature rule:
 	\begin{align*}
		Q\left(\mathbf{u}_h(\cdot,t);l\right) &= \sum_{i=0}^2 \omega_i \mathbf{f}(\mathbf{u}_h(x_j+\frac{h}{2}\xi_i,t)) \frac{\mathrm{d} \zeta_l(y)}{\mathrm{d}y}|_{y=\xi_i},
 	\end{align*}
 	with $\omega_0 = \omega_2 = 5/9$, $\omega_1 = 8/9$, and $\xi_0 = \sqrt{3/5}$, $\xi_1 = 0$, $\xi_2 = -\sqrt{3/5}$. Finally, system~\eqref{eq:DG_odesystem} gives us
 	\begin{equation*}
 		\begin{aligned}
 			\disp\frac{\mathrm{d}\mathbf{u}^{(0)}_j(t)}{\mathrm{d}t} &=
 			- \frac{1}{h} \left[ \Fnum(\mathbf{u}^-_{j+1/2}(t),\mathbf{u}^+_{j+1/2}(t)) - \Fnum(\mathbf{u}^-_{j-1/2}(t),\mathbf{u}^+_{j-1/2}(t)) \right],\\
 			\disp\frac{\mathrm{d}\mathbf{u}^{(1)}_j(t)}{\mathrm{d}t} &=
 			\frac{3}{h} \sum_{i=0}^2 \omega_i \mathbf{f}(\mathbf{u}_h(x_j+\frac{h}{2}\xi_i,t)) - \frac{3}{h} \left[ \Fnum(\mathbf{u}^-_{j+1/2}(t),\mathbf{u}^+_{j+1/2}(t)) + \Fnum(\mathbf{u}^-_{j-1/2}(t),\mathbf{u}^+_{j-1/2}(t)) \right],\\
 			\disp\frac{\mathrm{d}\mathbf{u}^{(2)}_j(t)}{\mathrm{d}t} &=
 			\frac{15}{h} \sum_{i=0}^2 \omega_i \xi_i \mathbf{f}(\mathbf{u}_h(x_j+\frac{h}{2}\xi_i,t)) - \frac{5}{h} \left[ \Fnum(\mathbf{u}^-_{j+1/2}(t),\mathbf{u}^+_{j+1/2}(t)) - \Fnum(\mathbf{u}^-_{j-1/2}(t),\mathbf{u}^+_{j-1/2}(t)) \right].
 		\end{aligned}
 	\end{equation*}
 	For linear systems $\mathbf{f}(\mathbf{u}) = A\mathbf{u}$, $A\in\mathbb{R}^{m\times m}$, the integral in~\eqref{eq:integral:discr} can be computed explicitly leading to
 	\begin{equation*}
 		\begin{aligned}
 			\disp\frac{\mathrm{d}\mathbf{u}^{(0)}_j(t)}{\mathrm{d}t} &=
 			- \frac{1}{h} \left[ \Fnum(\mathbf{u}^-_{j+1/2}(t),\mathbf{u}^+_{j+1/2}(t)) - \Fnum(\mathbf{u}^-_{j-1/2}(t),\mathbf{u}^+_{j-1/2}(t)) \right],\\
 			\disp\frac{\mathrm{d}\mathbf{u}^{(1)}_j(t)}{\mathrm{d}t} &=
 			\frac{6}{h} A \mathbf{u}_j^{(0)}(t) - \frac{3}{h} \left[ \Fnum(\mathbf{u}^-_{j+1/2}(t),\mathbf{u}^+_{j+1/2}(t)) + \Fnum(\mathbf{u}^-_{j-1/2}(t),\mathbf{u}^+_{j-1/2}(t)) \right],\\
 			\disp\frac{\mathrm{d}\mathbf{u}^{(2)}_j(t)}{\mathrm{d}t} &=
 			\frac{10}{h} A \mathbf{u}_j^{(1)}(t) - \frac{5}{h} \left[ \Fnum(\mathbf{u}^-_{j+1/2}(t),\mathbf{u}^+_{j+1/2}(t)) - \Fnum(\mathbf{u}^-_{j-1/2}(t),\mathbf{u}^+_{j-1/2}(t)) \right].
 		\end{aligned}
 	\end{equation*}
\end{example}

\subsection{Implicit time-integration}

The semi-discrete scheme~\eqref{eq:DG_odesystem} is an ODE system. Usually, one discretizes it using the total variation diminishing (TVD) Runge-Kutta (RK) time discretization~\cite{Shu1988}, which completes the definition of RKDG method in the explicit case. However, in this paper we are interested in formulating and analyzing implicit time approximations of the DG system. In particular, we focus on DIRK methods, with $s$ stages and general Butcher tableau
\begin{equation}\label{eq:dirk-tableau}
	\begin{array}{c|cccc}
		c_1 & a_{11} & 0 & \dots & 0 \\[1.5ex]
		c_2 & a_{21} & a_{22} & \dots & 0 \\[1.5ex]
		\vdots & \vdots & \vdots & \ddots & \\[1.5ex]
		c_s & a_{s1} & a_{s2} & \dots & a_{ss} \\[1ex]
		\hline
		&&&&\\[-1.8ex]
		& b_1 & b_2 & \dots & b_s
	\end{array}
\end{equation}
with consistency condition $\sum_{i=1}^s b_i = 1$. The other conditions on the coefficients of the tableau are chosen in order to have a desired order of accuracy.

The integration of~\eqref{eq:DG_odesystem} with a DIRK method of time-step $\Delta t$ leads to the following fully discrete numerical scheme:
\begin{equation}\label{eq:DG_DIRK}
	\begin{aligned}
	\mathbf{u}^{(l),n+1}_j &= \mathbf{u}^{(l),n}_j - (2l+1)\frac{\Delta t}{h} \sum_{i=1}^s b_i K_j^{(l),(i)} \\
	K_j^{(l),(i)} &= - Q\left(\mathbf{u}_h(\cdot,t^{(i)});l\right) + \left[ \Fnum(\mathbf{u}^{-,(i)}_{j+1/2},\mathbf{u}^{+,(i)}_{j+1/2}) - (-1)^l \Fnum(\mathbf{u}^{-,(i)}_{j-1/2},\mathbf{u}^{+,(i)}_{j-1/2}) \right],
	\end{aligned}
\end{equation}
where $\mathbf{u}^{\pm,(i)}_{j+1/2}$ are the left and right limits of the discontinuous solution $\mathbf{u}_h$ at the cell interface $x_{j+1/2}$ and time $t^{(i)} = (n + c_i) \Delta t$, whose moments are
\begin{equation}\label{eq:stage_values}
	\mathbf{u}^{(l),(i)}_j = \mathbf{u}^{(l),n}_j - (2l+1)\frac{\Delta t}{h} \sum_{\kappa=1}^i a_{i\kappa} K_j^{\Rtwo{(l)},(\kappa)}.
\end{equation}
Here and in the following, in order not to burden the notation we use a uniform time-step $\Delta t$. Nevertheless, this choice is not restrictive since the scheme can be formulated for a non-uniform time-step.

The advantage of DIRK schemes is that the implicit computation of a given stage value~\eqref{eq:stage_values} can be performed sequentially from $i=1$ to $i=s$. Therefore, at each time-step one has to solve $s$ systems of $m(p+1)N$ nonlinear equations:
\begin{equation}\label{eq:DIRK:stage}
	\mathbf{u}_j^{(l),(i)} 
	+ (2l+1)\frac{a_{ii}\Delta t}{h} K_j^{\Rtwo{(l)},(i)} 
	- \mathbf{u}_j^{(l),n} 
	+ (2l+1)\frac{\Delta t}{h} \sum_{\kappa=1}^{i-1} a_{ki} K_j^{\Rtwo{(l)},(\kappa)} = \mathbf{0}, \quad j=1,\dots,N, \ l=0,\dots,p.
\end{equation}
Above we have highlighted the term $K^{(i)}$, which makes the system for the $i$-th stage value nonlinear. In fact, the other space approximation terms $K^{(\kappa)}$, $\kappa=1,\dots,i-1$, are already available, thanks to the structure of DIRK schemes. Observe that the stencil of~\eqref{eq:DIRK:stage} is given by the cells $j-1$, $j$ and $j+1$, independently of the order of the DG approximation. Moreover, up to now, the non-linearity of~\eqref{eq:DIRK:stage} is only due to the nonlinear flux function $\mathbf{f}$, which enters in the computation of the numerical fluxes and of the quadrature rule.

\section{Dissipation-dispersion analysis of DIRK-DG schemes} \label{sec:analysis}

We analyze dissipation and dispersion errors within fully-discrete DIRK-DG formulations. The study will provide us knowledge on the best combination between space and time approximation.

To this end, we consider the scalar ($m=1$) linear advection equation in which $f=a u$, with $a$ scalar constant representing the wave speed. Specifically, we consider the initial value problem
\begin{equation}\label{eq:adv_lin}
	\begin{cases}
		u_t + a \, u_x = 0, & t>0 \\
		u(x,0) = e^{i\kk x},
	\end{cases}
\end{equation}
with periodic boundary conditions. The exact solution is
\begin{equation}\label{eq:adv_lin_sol}
	u(x,t) = e^{i\kk(x- a t)}
\end{equation}
which represents a sinusoidal wave train with a wave number $\kk$ and a frequency $a\kk$. In order to compare numerical and exact solutions, we consider the \Rtwo{$L^2$-projection} $\Pi_{V_{h}^p}$ of~\eqref{eq:adv_lin_sol} to the DG space $V_{h}^p$, so that we can write
\begin{equation}\label{eq:sol_ex_proj}
\Pi_{V_{h}^p}(e^{i\kk(x- a t)}) = \sum_{l=0}^p c_l(\kk h) e^{i\kk(x_j-at)}\, v^{(l)}_j(x) \quad \mbox{for } x\in I_j,
\end{equation}
with
\begin{equation}\label{eq:cl}
c_l(\kk h) = \disp\frac{1}{\int_{-1}^1 \zeta_l(\xi)^2 \mathrm{d}\xi}\int_{-1}^1 e^{i\kk\xi h/2}\zeta_l(\xi)\mathrm{d}\xi.
\end{equation}

The scheme propagates the projection of the initial condition on the DG space. Thus, in each cell we will study how the projection $\Pi_{V_{h}^p}(e^{i\kk x})$ is updated. The numerical solution defined on the DG space $V_h^p$ has the form given in \eqref{eq:uh}, with $u^{(l)}_j(t)$ degrees of freedom of the DG approximation for $l = 0,\cdots, p$.
We seek a numerical solution that closely approximates the projection~\eqref{eq:sol_ex_proj} of the exact solution. To achieve this, we use an Ansatz to express the degrees of freedom in the form
\begin{equation}\label{eq:sin_sol_num}
u^{(l)}_j(t) = c_l(\kk h) e^{i(\kk x_j-\omega t)}, \quad l=0,\ldots,p,
\end{equation}
where $\omega\in\C$ is a function of the wave number $\kk$ and defines the \emph{numerical frequency} of the scheme. We should expect that the numerical frequency matches $\kk a$ to recover the projection of the exact solution. Naively, one could expect that all modes contribute in the same way to approximate $\kk a$ $\forall \, l=0,\dots,p$. We will see, however, that the situation is far more complicated.

The distance between $\omega$ and $\kk a$ provides the error introduced by the numerical scheme. By splitting $\omega$ into real and imaginary parts, $\Real(\omega)$ and $\Imag(\omega)$, respectively, we can distinguish between dispersion and dissipation errors. In fact, we can write the Ansatz~\eqref{eq:sin_sol_num} as
$$
u^{(l)}_j(t) = c_l(\kk h) e^{i(\kk x_j-\omega t)} = c_l(\kk h) e^{\Imag(\omega)t}\, e^{i(\kk x_j - \Real(\omega)t)}, \quad l=0,\ldots,p.
$$
Thus, $\Real(\omega)$ represents the numerical speed, which characterizes the dispersion error, and to be exact we should have $\Real(\omega)= \kk a$. Instead, $\Imag(\omega)$ characterizes the dissipation error, and we would have $\Imag(\omega)=0$. 
When $\Real(\omega)\ne a\kk$ we say that the scheme is dispersive, whereas $\Imag(\omega)<0$ gives a dissipative solution. By taking into account how fine the mesh is with respect to the wave number $\kk$, we further define $K = \kk h$ and $\Omega = \omega h$. Then, the exact dispersion and dissipation relation can be expressed as $\mathfrak{Re}(\Omega)=aK$ and $\mathfrak{Im}(\Omega)=0$, respectively.

In the following, we aim to estimate the dispersion error $|\mathfrak{Re}(\Omega)-aK|$ and the dissipation error $\mathfrak{Im}(\Omega)$ of DIRK-DG methods with respect to the Courant number $C=|a|\frac{\dt}{h}$ for the wave numbers that can be resolved on the grid, which means $K\ll 1$.

\subsection{Numerical dissipation and dispersion of fully-discrete DG schemes} \label{sec:analysis:measures}

First, we give the recipe for computing the numerical frequency of a full space-time discretization of the linear equation~\eqref{eq:adv_lin} with a scheme of generic order $p+1$. Then, we specify the results for DIRK-DG methods for $p=1,2$.

Let us consider a generic numerical flux function
\begin{equation*}
	\Fnum = \Fnum(u^-,u^+) = a\left(b^+u^-+b^-u^+\right), \quad b^+=\frac{a+\beta|a|}{2a}, \ b^-=\frac{a-\beta|a|}{2a},
\end{equation*}
parameterized by $\beta\in[0,1]$. In particular, for the linear flux function $f(u)=au$, we use $\beta=1$ which recovers the Upwind numerical flux.
Then, the scheme~\eqref{eq:DG_odesystem} becomes
\begin{equation} \label{eq:scheme:linEq}
	\begin{array}{c}
		\disp\frac{\dx}{2}Z_l\, \frac{\mathrm{d}u_j^{(l)}(t)}{\mathrm{d}t} - a \sum_{m=0}^p Z^\prime_{lm}\ u^{(m)}_j(t)  
		+ a\,\sum_{m=0}^p \Big[ \left(b^+ u^{(m)}_j(t) + (-1)^m b^- u^{(m)}_{j+1}(t) \right)
		\medskip\\
		- (-1)^l \left(b^+ u^{(m)}_{j-1}(t) + (-1)^m b^- u^{(m)}_{j}(t) \right) \Big] = 0, \quad l = 0,\ldots, p,
	\end{array}
\end{equation}
with $Z_l = \int_{-1}^1 \zeta_l(\xi)^2\mathrm{d}\xi = 2/(2l+1)$, $Z^\prime_{lm} = \int_{-1}^1 \zeta_m(\xi)\zeta_l^\prime(\xi)\mathrm{d}\xi$.

Let $\mathbf{U}_j(t) = (u^{(0)}_j(t),\ldots,u^{(p)}_j(t))^T$ be the vector of the unknowns, i.e.~the degrees of freedom of the DG approximation in space. System~\eqref{eq:scheme:linEq} can be written in compact form as the following semi-discrete system of ODEs:
\begin{equation}\label{eq:scheme_lin_compact}
\frac{\mathrm{d}\mathbf{U}_j(t)}{\mathrm{d}t} = \frac{2\,a}{\dx}\left[D^-\,\mathbf{U}_{j-1}(t)+D\,\mathbf{U}_j(t)+D^+\,\mathbf{U}_{j+1}(t)\right],
\end{equation}
where $D^\pm$, $D$ are $p\times p$ time-independent matrices whose elements are
\begin{equation}
D^-_{lm} = \frac{b^+ (-1)^l}{Z_l}, \quad D_{lm} = \frac{Z^\prime_{lm}-b^+  + b^- (-1)^{m+l}}{Z_l}, \quad D^+_{lm} = \frac{-b^- (-1)^m}{Z_l}.
\end{equation}
Note that these matrices do not depend on $\kk$. We consider~\eqref{eq:scheme_lin_compact} endowed with the initial condition $\mathbf{U}_j(0)$ obtained by the projection of the initial datum $u(x,0) = e^{i\kk x}$ to the DG solution space $V_h^p$. Therefore, from~\eqref{eq:sol_ex_proj} we obtain $\mathbf{U}_j(0) = \boldsymbol{\mu}(\kk h) e^{i\kk x_j}$ where $\boldsymbol{\mu}(\kk h)=\left(c_0(\kk h),\ldots,c_p(\kk h)\right)^T$ with $c_l(\kk h)$, $l=0,\ldots,p$, defined in~\eqref{eq:cl}. In the following, we prefer not to burden the notation, and thus we omit the dependence of $\boldsymbol{\mu}$ on $\kk h$.

We study how~\eqref{eq:scheme_lin_compact} propagates the projection of the initial datum to get the degrees of freedom of the DG approximation. We recall that, according to Ansatz~\eqref{eq:sin_sol_num}, we seek solutions of the following form:
\begin{equation}\label{eq:sol_num_G}
\mathbf{U}_j(t) = \mathbf{W}(t) e^{i\kk x_j}, \quad t>0,
\end{equation}
and we expect that $\mathbf{W}(t)$ provides an accurate approximation of $\boldsymbol{\mu} e^{-i\kk a t}$. Substituting~\eqref{eq:sol_num_G} in~\eqref{eq:scheme_lin_compact}, we obtain the following system of ODEs for the time-dependent unknowns $\mathbf{W}(t)$:
\begin{equation}\label{eq:semi-discr-sys-g}
\disp\frac{\mathrm{d}\mathbf{W}(t)}{\mathrm{d}t} = \frac{2\,a}{\dx}\left(D^-\,e^{-i\kk\dx}+D\,+D^+\,e^{i\kk\dx}\right)\mathbf{W}(t),
\end{equation}
endowed with the initial condition $\mathbf{W}(0) = \boldsymbol{\mu}$.

\begin{remark}
	The amplification matrix $A = \frac{2\,a}{\dx}\left(D^-\,e^{-i\kk\dx}+D\,+D^+\,e^{i\kk\dx}\right)$, that defines the right-hand side of the semi-discrete DG system~\eqref{eq:semi-discr-sys-g}, has been already investigated in the literature. We refer, e.g., to~\cite{Guo2013} where it is shown that the expansion of the eigenvalues of $A$ as $h\to0$ with $a=1$ is:
	\begin{description}
		\item[$\mathbb{P}_1$ DG case:] $e_1 = -i\kk - \frac{\kk^4}{72} h^3 - \frac{i\kk^5}{270} h^4 + O(h^5)$ and $e_2 = -\frac{6}{h} + O(1)$;
		\item[$\mathbb{P}_2$ DG case:] $e_1 = -i\kk - \frac{\kk^6}{7200} h^5 - \frac{i\kk^7}{42000} h^6 + O(h^7)$ and $e_{2,3} = \frac{-3\pm\sqrt{51}i}{h} + O(1)$.
	\end{description}
	This proves that $A$ is diagonalizable with $p+1$, $p=1,2$, distinct eigenvalues. One of these eigenvalues is physically relevant since it approximates $-i\kk$ with dissipation error of order $2p+1$ and dispersion error of order $2p + 2$.
\end{remark}

Since the system ODEs is linear, applying a one-step time-discretization method of time-step $\Delta t > 0$ we write the solution of~\eqref{eq:semi-discr-sys-g} as
\begin{equation}
\mathbf{W}^n = M\,\mathbf{W}^{n-1}=\ldots=M^n\,\mathbf{W}(0), \quad n>0,
\end{equation} 
where $\mathbf{W}^n \approx \mathbf{W}(t^n)$, $t^n = n\Delta t$, and $M$ is a $(p+1)\times(p+1)$ matrix representing the fully-discrete numerical scheme. Therefore, the numerical solution of~\eqref{eq:semi-discr-sys-g} can be written in terms of the initial condition as $\mathbf{W}^n = M^n \boldsymbol{\mu}$. Note that the matrix $M$ depends on $\kk$, and hence its eigenvalues and eigenvectors will also depend on $\kk$. Moreover, assuming $M$ is diagonalizable, we can express the coefficients $\boldsymbol{\mu}$ with respect to the basis of the normalized eigenvectors $\{\boldsymbol{\theta}_l\}_{l=0}^p$ of $M$, i.e. 
\begin{equation}\label{eq:coefCm}
\sum_{l=0}^p \nu_l\, \boldsymbol{\theta_{l}} = \boldsymbol{\mu},
\end{equation}
and then $\mathbf{W}^n$ can be formulated as 
\begin{equation}
\mathbf{W}^n = M^n \boldsymbol{\mu} = \sum_{l=0}^p \nu_l\,M^n\,\boldsymbol{\theta}_l = \sum_{l=0}^p \nu_l\,\lambda^n_l\,\boldsymbol{\theta}_l,
\end{equation}
which is a superposition of eigenmodes along the eigenvectors of $M$ weighted by the expansion coefficients \Rone{$\{\nu_l\}_{l=0}^p$ of the initial condition in \eqref{eq:coefCm}}, and where $\{\lambda_l\}_{l=0}^p$ are the eigenvalues of $M$. Recalling Ansatz~\eqref{eq:sol_num_G}, we have that the numerical solution of the degrees of freedom $\mathbf{U}_j$ at time $t^n$ is
\begin{equation}
\mathbf{U}_j^n = \sum_{l=0}^p \nu_l\,\lambda^n_l\,\boldsymbol{\theta}_l\,e^{i\kk x_j} \quad\Rightarrow\quad u^{(m),n}_j = \sum_{l=0}^p \nu_l\,\lambda^n_l \theta_{ml}\,e^{i\kk x_j} \quad m=0,\ldots,p,
\end{equation}
and, finally, on each element $I_j$ the DG polynomial writes as 
\begin{equation}\label{eq:uh_sin}
u_h(x,t^n) = \sum_{m=0}^p \left( \sum_{l=0}^p \nu_l\,\lambda^n_l \theta_{ml}\,e^{i\kk x_j}\right) v^{(m)}_j(x) = \left\langle \mathbf{U}_j^n,\mathbf{v}_j(x) \right\rangle = \left\langle \mathbf{W}^n,\mathbf{v}_j(x) \right\rangle e^{i\kk x_j},
\end{equation}
where $\mathbf{v}_j(x) = (v^{(0)}_j(x),\dots,v^{(p)}_j(x))^T$ and $\langle \cdot,\cdot \rangle$ is the dot-product. Equation~\eqref{eq:uh_sin} shows that $\mathbf{U}_j^n = \mathbf{W}^n e^{i\kk x_j}$ is indeed a solution of the DG system, with $\mathbf{W}^n$ not depending on the cell.

In order, to compare the DG polynomial on the element $I_j$ with the projection on the DG space of the exact solution, we write also~\eqref{eq:sol_ex_proj} using~\eqref{eq:coefCm} to get
\begin{align}\label{eq:u_sin_ex}
\Pi_{V_{h}^p}(e^{i\kk(x- a t)}) = \sum_{m=0}^p \left(\sum_{l=0}^p \nu_l \theta_{ml}e^{i\kk(x_j-at^n)}\right) v^{(m)}_j(x) = \left\langle \sum_{l=0}^p \nu_l e^{-i\kk at^n} \boldsymbol{\theta}_{l},\mathbf{v}_j(x) \right\rangle e^{i\kk x_j}.
\end{align}
The relations for dissipation and dispersion introduced by the numerical scheme are determined from the absolute error between~\eqref{eq:uh_sin} and~\eqref{eq:u_sin_ex}.
\Rone{
\begin{prop} Let $M$ be the $(p+1)\times(p+1)$ matrix representing the fully-discrete numerical scheme applied to the linear problem \eqref{eq:adv_lin}, with exact solution \eqref{eq:adv_lin_sol}.
Let $M$ be a diagonalizable matrix, and let $\{(\lambda_l,\boldsymbol{\theta}_l)\}_{l=0}^p$ denote the set of eigenvalues of $M$ and their corresponding normalized, linearly independent eigenvectors. At a fixed time $t^n=n\dt$, the \emph{error} between the approximate solution \eqref{eq:uh} and the exact solution, is defined by the following components,
\begin{equation}\label{eq:err}
	\mathcal{E}^n_m = \sum_{l=0}^p \nu_l\, \boldsymbol{\theta}_{ml} \left( \lambda_l^n - e^{-i\kk at^n} \right), \, m=0,\ldots,p,
\end{equation}
where the coefficients $\{\nu_l\}_{l=0}^p$ are defined in \eqref{eq:coefCm} with $\boldsymbol{\mu}(\kk h)=\left(c_0(\kk h),\ldots,c_p(\kk h)\right)^T$, $c_l(\kk h)$, $l=0,\ldots,p$, given in~\eqref{eq:cl}.
\end{prop}
\begin{proof}
Since $M$ is diagonalizable, the approximate solution \eqref{eq:uh} can be expressed as in \eqref{eq:uh_sin}. Then,
\eqref{eq:err} is obtained by subtracting
\eqref{eq:uh_sin} and \eqref{eq:u_sin_ex}.
\end{proof}
Espression \eqref{eq:err} shows that each eigenmode imparts a specific level of dissipation and dispersion to the numerical solution.} Notice that in a single time-step the numerical scheme neither dissipates nor disperses if either $\lambda_l = e^{-i\kk a \Delta t}$ \ $\forall\,l=0,\dots,p$, or there exists $\hat{l}\in\{0,\dots,p\}$ such that $\lambda_{\hat{l}} = e^{-i\kk a \Delta t}$ and $\nu_l = 0$ \ $\forall\,l\neq\hat{l}$.

Let us write $\lambda_l = e^{-i\Omega_l\Delta t}$, with $\Omega_l = \omega_l h$ and where $\omega_l$ defines the numerical frequency. We can write the real and imaginary parts of $\Omega_l$ in terms of the eigenvalues of $M$ as
\begin{equation}\label{eq:full-discr-omega}
\Imag(\Omega_l) = \frac{h}{2\dt}\ln\left(\mathfrak{Re}(\lambda_l)^2+\Imag(\lambda_l)^2\right), \quad \mathfrak{Re}(\Omega_l) = -\frac{h}{\dt}\arctan(\frac{\Imag(\lambda_l)}{\mathfrak{Re}(\lambda_l)}).
\end{equation}
Thus, we observe that the numerical dissipation in an eigenmode is characterized by $\Imag(\Omega_l) < 0$, whereas numerical dispersion is present if $\mathfrak{Re}(\Omega_l) - aK \neq 0$. To measure how much a single eigenmode contributes to the numerical dissipation and dispersion, we observe that, \Rone{for each $m=0,\ldots,p$
\begin{equation}\label{eq:norm:err}
	\norm{\mathcal{E}^n_m}^2 \leq \sum_{l=0}^p |\nu_l|^2 | e^{-i(\Omega_l-aK)n\Delta t} - 1|^2,
\end{equation}
}
and thus, similarly to~\cite{Asthana2015}, we consider the \textit{relative energy} among modes
\begin{equation}\label{eq:energy}
	\beta_l = \frac{|\nu_l|^2}{\sum_{m=0}^p |\nu_m|^2}, \quad l=0,\ldots,p,
\end{equation}
which allows us to introduce the following measures of dissipation and dispersion:
\begin{description}
\item[Measure of dispersion]
\begin{equation}\label{eq:Mdisp}
	\mathcal{M}_{\text{disp}} = \disp\frac{1}{p+1}\sum_{l=0}^p \beta_l |\mathfrak{Re}(\Omega_l)-aK|
\end{equation}
\item[Measure of dissipation]
\begin{equation}\label{eq:Mdiss}
	\mathcal{M}_{\text{diss}}=\disp\frac{1}{p+1}\sum_{l=0}^p \beta_l \mathfrak{Im}(\Omega_l).
\end{equation}
\end{description}
We recall the dependence of $\beta_l$ on $\kk h$. In fact, by relation~\eqref{eq:coefCm}, $\beta_l$ depends on the initial condition which is expanded in the eigenvector basis to obtain the initial distribution of energy among the numerical modes as a function of the wave number. Quantities~\eqref{eq:Mdisp} and~\eqref{eq:Mdiss} are weighted by $1/(p+1)$ in order to make a fair comparison between schemes of different order.

\subsubsection{The case of DIRK-DG approximations}
The previous analysis holds true for a general one-step time-integration method. Let us now consider $s$-stage DIRK schemes~\eqref{eq:dirk-tableau} for the time-integration of the semi-discrete system of ODEs~\eqref{eq:semi-discr-sys-g} and let us denote $A = \frac{2\,a}{\dx}\left(D^-\,e^{-i\kk\dx}+D\,+D^+\,e^{i\kk\dx}\right)$ the matrix of the space approximation in~\eqref{eq:semi-discr-sys-g}. Then, the fully-discrete scheme reads
\begin{equation}\label{eq:dirk-scheme}
\begin{aligned}
\mathbf{W}^{n+1} &= \mathbf{W}^n + \dt \disp\sum_{r=1}^s b_r A\,\mathbf{Y^r}, \quad n\geq 0 \\
\mathbf{Y^r} &= \mathbf{W}^n + \dt \disp\sum_{q=1}^r a_{rq} A\,\mathbf{Y^q}, \quad r=1,2,\ldots,s.
\end{aligned}  
\end{equation}
Since $a_{rq}=0$ for all $q>r$, we then have
$$
\left(\Id-\dt\,a_{rr}\,A\right)\mathbf{Y^r} = \mathbf{W}^n + \dt\disp\sum_{q=1}^{r-1} a_{rq} A\,\mathbf{Y^q},
$$
where $\Id$ is the $(p+1)\times(p+1)$ identity matrix. Denoting $B_r = \Id-\dt\,a_{rr}\,A$, we have
$$
\mathbf{Y^r} = B_r^{-1}\left(\mathbf{W}^n + \dt\disp\sum_{q=1}^{r-1} a_{rq} A\,\mathbf{Y^q}\right) = B_r^{-1}\mathbf{W}^n + \dt\disp\sum_{q=1}^{r-1} a_{rq} B_r^{-1}\,A\,\mathbf{Y^q}. 
$$
For DIRK methods with $s\leq 4$ stages, the general formulation of $M$ such that $\mathbf{W}^{n+1} = M \mathbf{W}^{n}$ is given by
\begin{equation} \label{eq:matrix:DIRKDG}
\begin{aligned}
	M &= \Id + \Delta t \sum_{i=1}^s b_i A B_i^{-1} + \Delta t^2 \sum_{i=2}^s \sum_{j=1}^{i-1} b_i a_{ij} A B_i^{-1} A B_j^{-1} + \Delta t^3 \sum_{i=3}^s \sum_{j=2}^{i-1} \sum_{k=1}^{j-1} b_i a_{ij} a_{jk} A B_i^{-1} A B_j^{-1} A B_k^{-1} \\
	&+ \Delta t^4 b_s a_{s,s-1} a_{s-1,s-2} a_{s-2,s-3} A B_s^{-1} A B_{s-1}^{-1} A B_{s-2}^{-1} A B_{s-3}^{-1}.
\end{aligned}
\end{equation}

Let us give an example for the $\mathbb{P}_1$ DG approximation coupled with a $2$-stage second-order DIRK scheme and the upwind numerical flux. In this case we have $\mathbf{W}^{n+1} = M \mathbf{W}^n$ with
$$
	M = \Id + \dt A (b_1 B_1^{-1}+ b_2 B_2^{-1})+ \dt^2 b_2 a_{21}A B_2^{-1} A B_1^{-1}
$$
	where $\Id$ is the $2\times 2$ identity matrix and
$$
	A = -\frac{a}{h}\begin{bmatrix} 1-e^{-i\kk h} & 1-e^{-i\kk h} \\[1ex] -3+3e^{-i\kk h} & 3+3e^{-i\kk h} \end{bmatrix}
$$
is the matrix of the space approximation. If we assume that the DIRK scheme is characterized by the Butcher tableau given in Table~\ref{tab:dirk22} in Section~\ref{sec:analysis:numerics}, then the matrices $B_1$ and $B_2$ are the same for all stages, and
\begin{align*}
	M &= \Id + \dt A B^{-1} + \frac{1-2\gamma}{2} \dt^2 A B^{-1} A B^{-1}, \\
	B &:= B_1 = B_2 = \Id - \gamma \Delta t A.
\end{align*}
The amplification matrix $M$ is diagonalizable with two distinct eigenvalues. Indeed, taking $a=1$, we get
\begin{align*}
\lambda_0 &= 1 - i \kk \Delta t - \frac{\kk^2}{2} \Delta t^2 + i\kk^3 (\gamma-\gamma^2) \Delta t^3 - \frac{\kk^4}{72}\Delta t h^3 + \frac{3\kk^4}{2} (3\gamma^2-2\gamma^3) \Delta t^4 + O(\Delta t^5),\\
\lambda_1 &= \frac{h^3 - h^2 \Delta t (6-18\gamma) + h \Delta t^2 (18 - 108 \gamma + 108\gamma^2) + \Delta t^3 (108 \gamma - 432 \gamma^2 + 216 \gamma^3) + 3i \kk h^3 \Delta t - \kk (18i-54i\gamma) h^2 \Delta t^2}{(6 \Delta t \gamma + h)^3} \\ &+ O(\Delta t^2).
\end{align*}
We observe that $\lambda_0$ is the \emph{principal} eigenvalue since it approximates $e^{-i\kk\Delta t}$ with dispersion error of order $4$ and dissipation error of order $3$. Instead, the second eigenvalue $\lambda_1$ introduces spurious numerical effects. To investigate the relative energy~\eqref{eq:energy}, we first observe that the expansion of the coefficients $c_l$, $l=0,1$, of the initial condition is
$$
	c_0 = \frac{2\sin(\frac{h\kk}{2})}{h\kk} = 1-\frac{h^2\kk^2}{24}+O(h^4), \quad c_1 = -\frac{6i(h\kk\cos(\frac{h\kk}{2}) - 2\sin(\frac{h\kk}{2}))}{h^2\kk^2} = \frac{ih\kk}{2} + O(h^3).
$$
Then, their expansion with respect to the basis of the normalized eigenvectors becomes
\begin{align*}
	\nu_0 = 1 + \frac{h^2\kk^2}{8} + \frac{i h^3\kk^3}{72} + O(h^4), \quad
	\nu_1 = -\frac{h^2\kk^2}{6} - \frac{ih^3\kk^3}{72} + O(h^4).
\end{align*}
Finally, we have the expansion of the relative energy as
\begin{align*}
	\beta_0 = 1 - \frac{h^4\kk^4}{36} - \frac{ih^5\kk^5}{216} + O(h^6), \quad
	\beta_1 = \frac{h^4\kk^4}{36} + \frac{ih^5\kk^5}{216} + O(h^6), 
\end{align*}
which shows that the principal eigenvalue is associated to the larger energy, and therefore it provides the main contribution to the dissipation-dispersion error of the scheme. Notice also that $\beta_0$ and $\beta_1$ do not depend on the free parameter of the DIRK scheme.


If we consider $\gamma=0.25$ the principal eigenvalue behaves as
\begin{align*}
\lambda_0 &= 1 - i\kk \Delta t - \frac{\kk^2}{2} \Delta t^2 + \frac{3i\kk^3}{16} \Delta t^3 - \frac{\kk^4}{72} \Delta t  h^3 + \frac{\kk^4}{16} \Delta t^4 + O(\Delta t^5),
\end{align*}
whereas, if $\gamma = 1-\sqrt{2}/2$, which corresponds to the L-stable DIRK scheme in~\cite{Pareschi2005}, we have
$$
\lambda_0 = 1 - i\kk \Delta t - \frac{\kk^2}{2} \Delta t^2 + \frac{204i\kk^3}{985} \Delta t^3 - \frac{\kk^4}{72} \Delta t  h^3 + \frac{306\kk^4}{985} \Delta t^4 + O(\Delta t^5).
$$

\Rone{
\begin{remark}
We can estimate the order of convergence of the scheme by using the error components defined in \eqref{eq:err} and the Taylor expansions mentioned above.
Using the upper bound 
$$
\|\mathcal{E}^n_m\|^2 \lesssim \sum_{l=0}^p |w_l|^2\,|\lambda_l^n - e^{-i \kk a t^n}|^2,
$$ 
we can estimate the order of convergence for each error term, $m=0,\ldots,p$.
On a single time step, for $n=1$ and from the expansions 
$\lambda_0 = 1 - i \kk \Delta t - \frac{\kk^2}{2} \Delta t^2 + O(\Delta t^3)$ and $\nu_0 = 1 + O(h^2)$, we get
$$
    \lambda_0 - e^{-i \kk a \Delta t} = O(\Delta t^3) \quad\mbox{and}\quad |\nu_0(\lambda_0 - e^{-i \kk a \Delta t})| = O(\Delta t^3).
$$
From $\lambda_1 = O\Big(\frac{\Delta t}{h}\Big) + O\Big(\frac{\Delta t^2}{h^2}\Big) + O\Big(\frac{\Delta t^3}{h^3}\Big) + \dots$ and $\nu_1 = O(h^2)$,
one obtains
$$
    |\nu_1 (\lambda_1 - e^{-i \kk a \Delta t})| \sim O(h^2).
$$
Therefore, with respect to the norm $\|\mathcal{E}^n_m\|$, $m=0,\ldots,p$ the dominant contributions are
$$
    \|\mathcal{E}_m\| \sim O(h^2) + O(\Delta t^3).
$$
Consequently, the method is of second order in space and third order in time. Furthermore: if $\Delta t \sim h$, the spatial error dominates, yielding overall order $2$; if $\Delta t \sim h^{2/3}$ so that the terms are balanced, the temporal order $\sim 3$ is achieved.
\end{remark}
}

For a $\mathbb{P}_2$ DG approximation, the matrix of the space approximation is given by
$$
	A = -\frac{a}{ h} \begin{bmatrix} 1-e^{-i\kk h} & 1-e^{-i\kk h} & 1-e^{-i\kk h} \\[1ex] -3+3e^{-i\kk h} & 3+3e^{-i\kk h} & 3+3e^{-i\kk h} \\[1ex] -5e^{-i\kk h}+5 & -5e^{-i\kk h}-5 & -5e^{-i\kk h}+5 \end{bmatrix}.
$$
Whereas, the amplification matrix, with $2\leq s \leq 4$, can be obtained from~\eqref{eq:matrix:DIRKDG}. However,
in this case the analysis of the eigenvalues is more complicated and therefore we proceed numerically. In fact, we can reason as before by computing the eigenvalues of $M$ and measuring the numerical dispersion and dissipation by means of~\eqref{eq:Mdisp} and~\eqref{eq:Mdiss}, respectively.

\subsection{Numerical investigation of dissipation and dispersion in DIRK-DG schemes} \label{sec:analysis:numerics}

Here, we investigate numerically the dissipation-dispersion errors of second and third-order DIRK-DG fully-discrete approximations. To this end, we numerically compute the eigenvalues of the corresponding space-time discretization matrices and evaluate the relations~\eqref{eq:Mdisp} and~\eqref{eq:Mdiss}. At the same time, we consider the linear advection equation~\eqref{eq:adv_lin} with $a=1$ on $x\in[-1,1]$ for regular and irregular initial data $u_0(x)$:
\begin{subequations}
	\begin{equation} \label{eq:init:sinhf}
		u_0(x) = -\sin(10\pi x),
	\end{equation}
	\begin{equation} \label{eq:init:smoothhf}
		u_0(x) = \sin(\pi x)+\sin(10\pi x) \exp(-20x^2)
	\end{equation}
	\begin{equation} \label{eq:init:doublestep}
		u_0(x) = \begin{cases} 0.5, & -0.25\leq x \leq 0.25, \\
			0, & \text{elsewhere}. \end{cases}
	\end{equation}
\end{subequations}
The objective is to illustrate how, while maintaining the same level of accuracy, the combination of DG in space and DIRK in time approximations yields vastly distinct outcomes for both regular and irregular initial data.

\Rtwo{Throughout the paper, we define
\begin{equation} \label{eq:dtCFL}
	\Delta t_{\mathrm{CFL}} := \frac{h}{(2p+1)\,\max_x \norm{\mathbf{f}'(\mathbf{u})}},
\end{equation}
where $\norm{\mathbf{f}'(\mathbf{u})}$ denotes the spectral radius of the Jacobian of the flux function. This corresponds to the classical CFL restriction for a discontinuous Galerkin scheme of degree $p$ with a Rusanov flux. In general, the precise stability bound of an explicit scheme depends also on the particular Runge-Kutta integrator employed. However, in this work we take \eqref{eq:dtCFL} as a reference definition of $\Delta t_{\mathrm{CFL}}$, independently of the specific Runge-Kutta scheme, in order to provide a consistent baseline for comparing explicit and implicit methods. We then introduce
\begin{equation} \label{eq:def:Courant}
	r := \frac{\Delta t_{\text{imp}}}{\Delta t_{\text{CFL}}},
\end{equation}
that is, the ratio between the implicit time-step $\Delta t_{\text{imp}}$ and the reference explicit stability limit $\Delta t_{\text{CFL}}$ running the same test with an explicit RKDG of matching order in space and time.
When $r=1$, implicit and explicit schemes are advanced with the same time-step.}

We compare the results obtained with one-parameter family of second-order DIRK methods with $s=2$ stages and the third-order DIRK methods with $s=3$ and $s=4$ stages described, respectively, by the Butcher tableaux in Table~\ref{tab:dirk22}, Table~\ref{tab:dirk33} and Table~\ref{tab:dirk43}.
\begin{table}[h!]
	\begin{minipage}{.25\linewidth}
		\centering
		\footnotesize
		\[\arraycolsep=5.0pt\def\arraystretch{1.4}
		\begin{array}{c|cc}
			\gamma & \gamma & 0 \\[1.5ex]
			1-\gamma & 1-2\gamma & \gamma \\[1.5ex]
			\hline
			& \tfrac12 & \tfrac12 
		\end{array}\]
		\caption{$2$-stage DIRK method. It is second-order accurate for any choice of $\gamma$.\label{tab:dirk22}}
	\end{minipage}
	\hfill
	\begin{minipage}{0.35\linewidth}
		\centering
		\footnotesize
		\[\arraycolsep=5.0pt\def\arraystretch{1.5}
		\begin{array}{c|ccc}
		\gamma & \gamma & 0 & 0 \\[1.4ex]
		1+\tfrac{\gamma}{2} & \tfrac{1-\gamma}{2} & \gamma & 0 \\[1.5ex]
		1 & -\tfrac32\gamma^2+4\gamma-\tfrac14 & \tfrac32\gamma^2-5\gamma+\tfrac54 & \gamma \\[1.5ex]
		\hline
		& -\tfrac32\gamma^2+4\gamma-\tfrac14 & \tfrac32\gamma^2-5\gamma+\tfrac54 & \gamma
		\end{array}\]
		\caption{$3$-stage DIRK method. It is third-order accurate for $\gamma=0.158983899988677$, $\gamma=0.435866521508459$, $\gamma=2.405149578502864$.\label{tab:dirk33}}
	\end{minipage}
	\hfill
	\begin{minipage}{0.25\linewidth}
		\centering
		\footnotesize
		\[\arraycolsep=5.0pt\def\arraystretch{1.5}
		\begin{array}{c|cccc}
			\tfrac12 & \tfrac12 & 0 & 0 & 0 \\[1.4ex]
			\tfrac23 & \tfrac16 & \tfrac12 & 0 & 0 \\[1.5ex]
			\tfrac12 & -\tfrac12 & \tfrac12 & \tfrac12 & 0 \\[1.5ex]
			1 & \tfrac32 & -\tfrac32 & \tfrac12 & \tfrac12 \\[1.5ex]
			\hline
			& \tfrac32 & -\tfrac32 & \tfrac12 & \tfrac12
		\end{array}\]
		\caption{$4$-stage third-order DIRK method.\label{tab:dirk43}}
	\end{minipage}
\end{table} 

\begin{figure}[t!]
	\centering
	\begin{subfigure}[b]{0.3\textwidth}
		\centering
		\includegraphics[width=\textwidth]{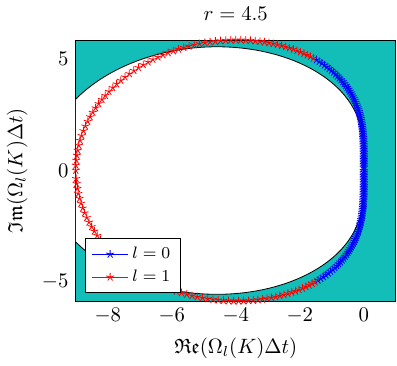}
		\includegraphics[width=\textwidth]{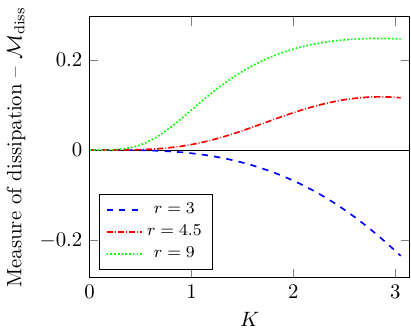}
		\caption{$2$-stage second-order DIRK method of Table~\ref{tab:dirk22} with $\gamma=0.2$.\label{fig:diss:DIRK2}}
	\end{subfigure}
	\hfill
	\begin{subfigure}[b]{0.3\textwidth}
		\centering
		\includegraphics[width=\textwidth]{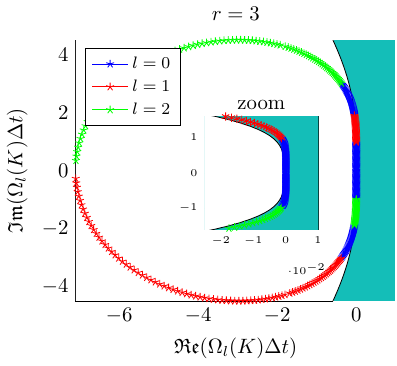}
		\includegraphics[width=\textwidth]{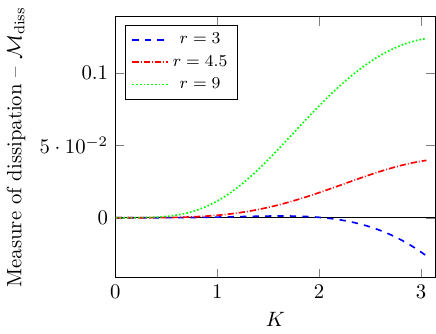}
		\caption{$3$-stage third-order DIRK method of Table~\ref{tab:dirk33} with $\gamma=0.158983899988677$.\label{fig:diss:DIRK3:gamma1}}
	\end{subfigure}
	\hfill
	\begin{subfigure}[b]{0.3\textwidth}
		\centering
		\includegraphics[width=\textwidth]{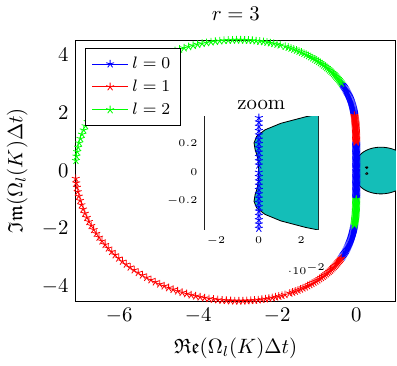}
		\includegraphics[width=\textwidth]{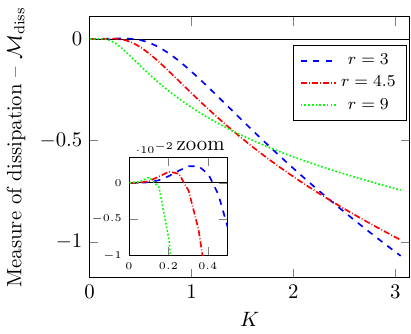}
		\caption{$3$-stage third-order DIRK method of Table~\ref{tab:dirk33} with $\gamma=2.405149578502864$.\label{fig:diss:DIRK3:gamma2}}
	\end{subfigure}
	\caption{We show the development of instability when the DIRK parameters do not yield A-stable schemes. Top row: The white area represents the stability regions of the DIRK methods of order $p+1=2,3$, whereas the $\star$ markers are obtained from the eigenvalues $\Omega_l$, $l=0,\dots,p$, \Rtwo{of the matrix coming from} the space $\mathbb{P}_p$ DG approximation as function of $K$. Bottom row: Measure of dissipation~\eqref{eq:Mdiss} computed for several Courant numbers.\label{fig:diss:unstable}}
\end{figure}
We observe that, for the $2$-stage DIRK in Table~\ref{tab:dirk22}, any choice of $\gamma$ satisfies second-order conditions. In particular, we consider $\gamma\geq 0.25$. As already observed, $\gamma = 1-\sqrt{2}/2 \approx 0.29289$ gives the L-stable DIRK method of Pareschi and Russo~\cite{Pareschi2005}. Values of $\gamma < 0.25$ do not guarantee A-stability of the method, as showed in Figure~\ref{fig:diss:DIRK2}. In the top panel, we plot the stability region of the $2$-stage DIRK in Table~\ref{tab:dirk22} with $\gamma=0.2$. The white area represents the stability region of the DIRK method which clearly does not contain the whole half-plane $\mathfrak{Re}(z)<0$, $z\in\mathbb{C}$. The $\star$ markers are the points $\mathcal{P}_l(K) = (\mathfrak{Re}(\Omega_l(K)\Delta t),\mathfrak{Im}(\Omega_l(K)\Delta t)$, $l=0,1$, in the complex plane, \Rtwo{where $\Omega_l(K)$ depends on the $l$-th eigenvalue of the matrix coming from the space $\mathbb{P}_1$ DG approximation obtained with the wave number $K=\kk h$.}
Here, $\Delta t_{\text{imp}}$ is computed with $r=4.5$, and we observe that for some $K$ the points $\mathcal{P}_l(K)$ lie outside the stability region of the DIRK method. In the bottom panel of Figure~\ref{fig:diss:DIRK2}, we show the measure of dissipation~\eqref{eq:Mdiss} for $r=3,4.5,9$. Notice that for small $r$ it is still possible to guarantee the stability of the scheme since the dissipation is negative. \Rone{We remark that the plots in the bottom panels of Figure~\ref{fig:diss:unstable} are restricted to $K \in [0,\pi]$ since the discrete Fourier spectrum is $2\pi$-periodic in $K$. Instead, in the top panels of Figure~\ref{fig:diss:unstable} we consider $K \in [0,2\pi]$ to show that the eigenvalue curves close, as typically expected.}

For the $3$-stage DIRK scheme in Table~\ref{tab:dirk33}, we observe that third-order conditions are satisfied for $\gamma=0.158983899988677$, $\gamma = 0.435866521508459$ and $\gamma=2.405149578502864$. The second value allows to recover the A-stable and stiffly accurate, and therefore L-stable, method of~\cite{Alexander1977}. We will consider only the value $\gamma=0.435866521508459$ since the other parameter values provide methods that are not A-stable. This is observed in the top panels of Figure~\ref{fig:diss:DIRK3:gamma1} and Figure~\ref{fig:diss:DIRK3:gamma2}, where we plot the stability region of the $3$-stage DIRK in Table~\ref{tab:dirk33} with $\gamma=0.158983899988677$ and $\gamma=2.405149578502864$, respectively. As in the second-order DIRK, the stability region is represented by the white area which does not contain the whole half-plane $\mathfrak{Re}(z)<0$, $z\in\mathbb{C}$. The $\star$ markers are the points $\mathcal{P}_l(K) = (\mathfrak{Re}(\Omega_l(K)\Delta t),\mathfrak{Im}(\Omega_l(K)\Delta t)$, $l=0,1,2$, in the complex plane, where $\Omega_l(K)$ is the $l$-th eigenvalue \Rtwo{of the matrix coming from} the space $\mathbb{P}_2$ DG approximation obtained with the wave number $K=\kk h$. Here, $\Delta t_{\text{imp}}$ is computed with $r=3$. The zoom shows that for some $K$ the points $\mathcal{P}_l(K)$ lie outside the stability regions of the DIRK methods. In the bottom panels of Figure~\ref{fig:diss:DIRK3:gamma1} and Figure~\ref{fig:diss:DIRK3:gamma2}, we show the measure of dissipation~\eqref{eq:Mdiss} for $r=3,4.5,9$. Notice that the instability of the schemes are highlighted also by the negative dissipation.

The L-stable third-order method of~\cite{Alexander1977}, obtained with the parameter $\gamma=0.435866521508459$, is compared in the following with the $4$-stage third-order accurate method given in Table~\ref{tab:dirk43} which is a Strong Stability Preserving (SSP) scheme given in~\cite{Santos2021}. 

We summarize in Table~\ref{tab:labels} the list of the schemes we will investigate in the following and in Section~\ref{sec:numerical:simulations}, with the corresponding labels. 

\begin{table}[h!]
	\caption{Labels and descriptions of the schemes tested}
	\begin{center}
		\begin{tabular}{c|c}
			\underline{Label} & \underline{Scheme description} 
			\\[2ex]
			$\DIRKdue$ & $\mathbb{P}_1$ DG coupled with the $2$-stage second-order DIRK in Table~\ref{tab:dirk22} with $\gamma=0.25$ 
			\\[2ex]
			$\DIRKlstab$ & $\mathbb{P}_1$ DG coupled with the $2$-stage second-order DIRK in Table~\ref{tab:dirk22} with $\gamma=1-\sqrt{2}/2$ 
			\\[2ex]
			$\DIRKlam{0.5}$ & $\mathbb{P}_1$ DG coupled with the $2$-stage second-order DIRK in Table~\ref{tab:dirk22} with $\gamma=0.5$ 
			\\[2ex]
			$\DIRKlstabTre$ & $\mathbb{P}_2$ DG coupled with the $3$-stage third-order DIRK in Table~\ref{tab:dirk33} with $\gamma=0.435866521508459$ 
			\\[2ex]
			$\DIRKQuattroTre$ & $\mathbb{P}_2$ DG coupled with the $4$-stage third-order DIRK in Table~\ref{tab:dirk43} 
			\\
		\end{tabular}
	\end{center}
	\label{tab:labels}
\end{table}

\begin{figure}[t!]
	\centering
	\begin{subfigure}[b]{0.48\textwidth}
		\centering
		\includegraphics[width=\textwidth]{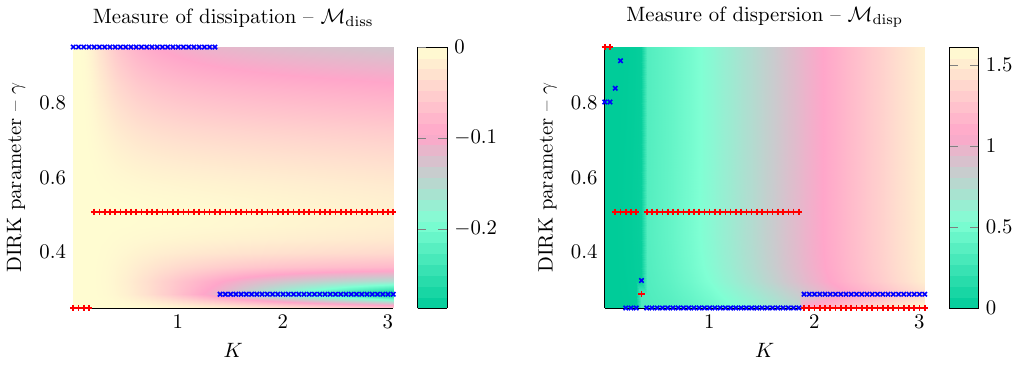}
		\includegraphics[width=\textwidth]{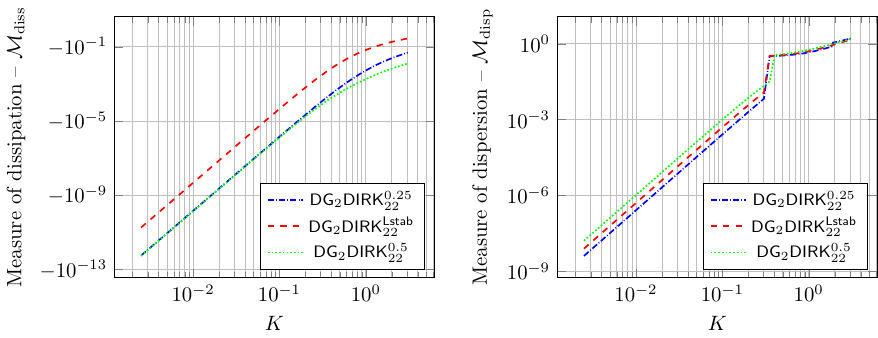}
		\caption{$r=15$.\label{fig:diss:disp:p1:cou5}}
	\end{subfigure}
	\hfill
	\begin{subfigure}[b]{0.48\textwidth}
		\centering
		\includegraphics[width=\textwidth]{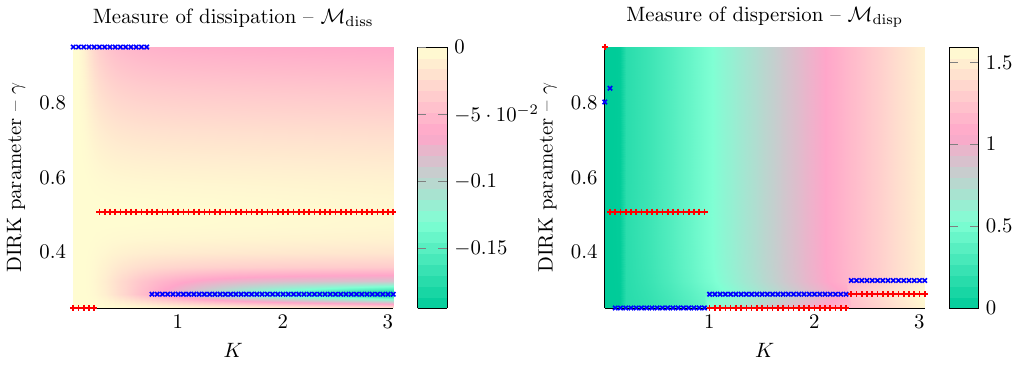}
		\includegraphics[width=\textwidth]{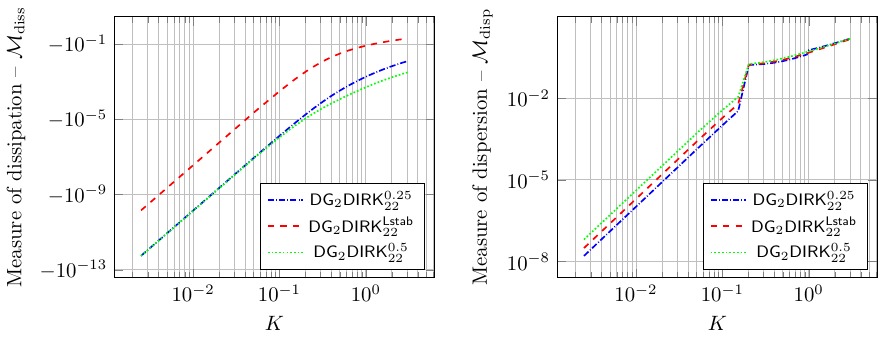}
		\caption{$r=30$.\label{fig:diss:disp:p1:cou10}}
	\end{subfigure}
	\caption{We explore the sensitivity of second-order schemes to the choice of the parameter $\gamma$. Top row: Measure of dissipation~\eqref{eq:Mdiss} and dispersion~\eqref{eq:Mdisp} for second-order DIRK-DG schemes as function of $K$ and $\gamma$. Blue $\times$ markers and red $+$ markers denote the values of the DIRK parameter for which the measure is minimum and maximum, respectively. Bottom row: Dissipation and dispersion measures obtained with DIRK parameters $\gamma=0.25,1-\sqrt{2}/2,0.5$.\label{fig:diss:disp:p1}}
\end{figure}

\begin{figure}[t!]
	\centering
	\begin{subfigure}[b]{0.48\textwidth}
		\centering
		\includegraphics[width=\textwidth]{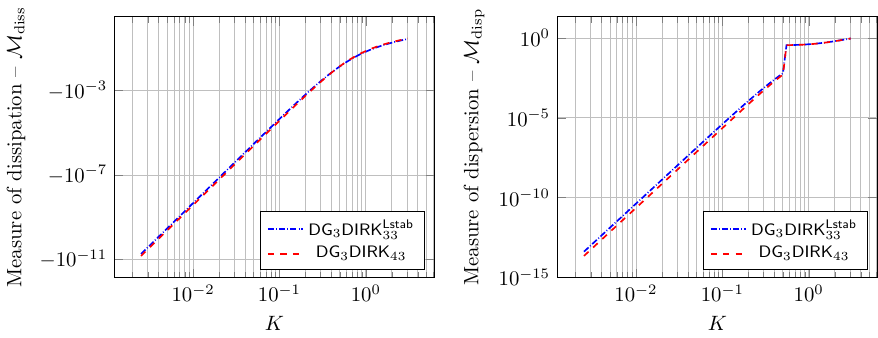}
		\caption{$r=15$.\label{fig:diss:disp:p2:cou3}}
	\end{subfigure}
	\hfill
	\begin{subfigure}[b]{0.48\textwidth}
		\centering
		\includegraphics[width=\textwidth]{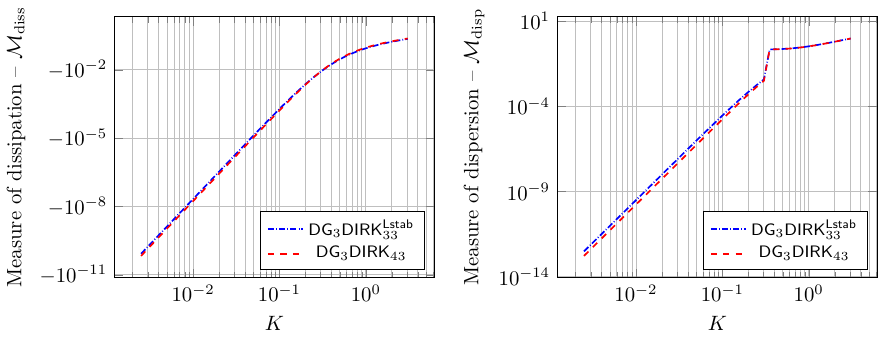}
		\caption{$r=25$.\label{fig:diss:disp:p2:cou5}}
	\end{subfigure}
	\hfill
	\begin{subfigure}[b]{0.48\textwidth}
		\centering
		\includegraphics[width=\textwidth]{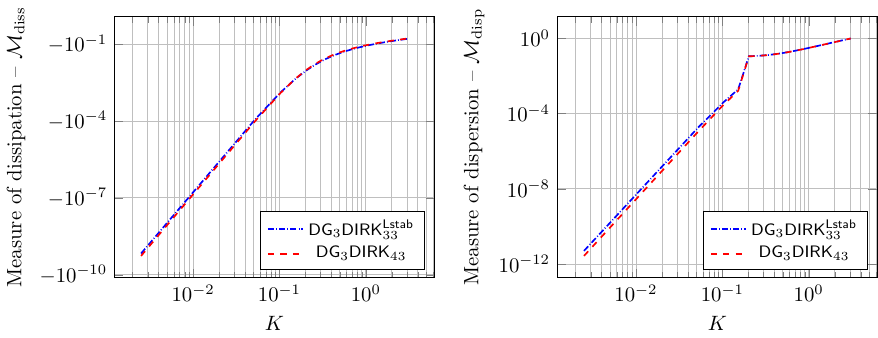}
		\caption{$r=50$.\label{fig:diss:disp:p2:cou10}}
	\end{subfigure}
	\caption{Measure of dissipation~\eqref{eq:Mdiss} and dispersion~\eqref{eq:Mdisp} for the two third-order DIRK-DG schemes given in Table~\ref{tab:dirk33} and Table~\ref{tab:dirk43}. Observe that the case $r=15$ in the top-left panel should be compared to the case $r=15$ of second-order schemes in Figure~\ref{fig:diss:disp:p1}, since the ratio between the time-step used and the time-step required for explicit stability is the same.\label{fig:diss:disp:p2}}
\end{figure}

Top row of Figure~\ref{fig:diss:disp:p1} shows the dissipation and dispersion errors as functions of the wave number $K = \kk  h$ with $ h=0.05$ and of the DIRK parameter $\gamma$ for $p=1$. The discretization of the parameter space is made with $20$ uniformly spaced points in $[0.25,0.95]$. The results are depicted for two different values of the ratio $r$, namely $r=15$ and $r=30$. We observe that the dissipation and dispersion errors are increasing with the wave number and the Courant number, as we expect. The DIRK parameter $\gamma=1-\sqrt{2}/2$ provides maximum diffusion error, whereas $\gamma=0.5$ has minimum diffusion but large dispersion error in the mid-range of the wave number, as observed also in the bottom row of Figure~\ref{fig:diss:disp:p1}.

In Figure~\ref{fig:diss:disp:p2}, we present the dissipation and dispersion measures of two selected DIRK schemes of order 3. Specifically, we examine the L-stable DIRK scheme from Table~\ref{tab:dirk33} with $\gamma=0.435866521508459$ and the SSP DIRK scheme from Table~\ref{tab:dirk43}. The results are displayed for three ratio numbers $r=15,25,50$.
We observe that the $4$-stage SSP DIRK scheme of order 3 is slightly less diffusive and less dispersive than the $3$-stage L-stable DIRK scheme of order 3. Additionally, both dissipation and dispersion measures increase with increasing Courant number.

Let us compare the results of second- and third-order schemes when they use the same ratio $r=15$ between the time-step used and the time-step required for explicit stability. We note that the measure of dispersion for third-order schemes is smaller than that for second-order schemes. Conversely, the measure of dissipation for third-order schemes is slightly larger than that for second-order schemes at small normalized wave numbers $K$. Given these observations, we expect minimal differences in the dissipation of solutions between DIRK-DG schemes with $p=1$ and $p=2$. However, we anticipate that second-order DIRK-DG schemes will produce more dispersive solutions, leading to poor approximation of the velocity of linear transport and the introduction of spurious oscillations in discontinuous profiles unless a limiter is introduced.

\begin{figure}[t!]
	\centering
	\begin{subfigure}[b]{\textwidth}
		\centering
		\includegraphics[width=\textwidth]{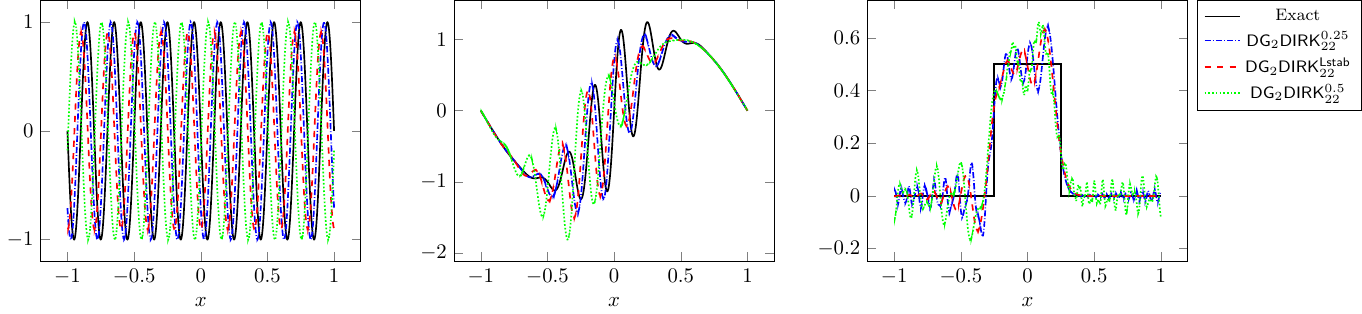}
		\caption{$r=15$.\label{fig:diss:disp:lintra:p1}}
	\end{subfigure}
	\begin{subfigure}[b]{\textwidth}
		\centering
		\includegraphics[width=\textwidth]{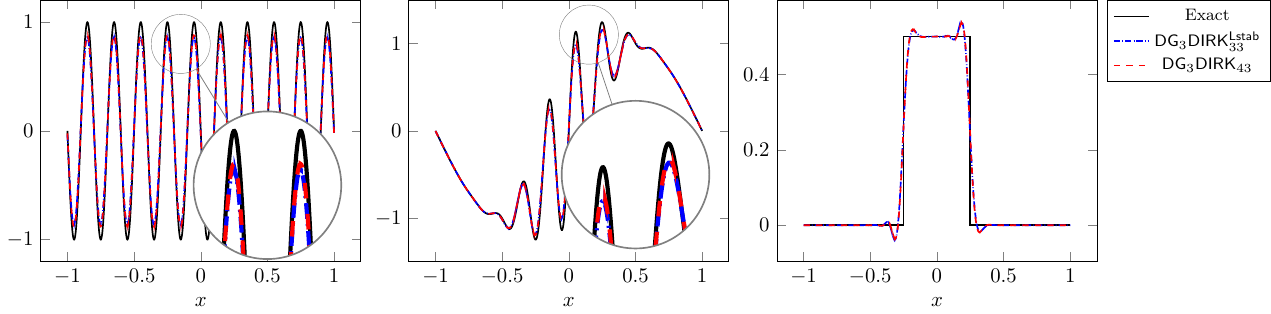}
		\caption{$r=15$.\label{fig:diss:disp:lintra:p2:cou3}}
	\end{subfigure}
	\caption{Solutions to the linear advection equation with speed $a=1$, final time $t=2$ and initial conditions~\eqref{eq:init:sinhf} (left column), \eqref{eq:init:smoothhf} (middle column), and~\eqref{eq:init:doublestep} (right column). All the solutions are obtained with $N=400$ cells.\label{fig:diss:disp:lintra}}
\end{figure}

We compare all these schemes on the solution of the linear advection equation~\eqref{eq:adv_lin} with speed $a=1$ and initial conditions~\eqref{eq:init:sinhf}, \eqref{eq:init:smoothhf} and~\eqref{eq:init:doublestep}. See Figure~\ref{fig:diss:disp:lintra}. Here, $r=15$ both for second- and third-order schemes, so that we can compare them when they have the same ratio between the time-step used and the time-step for the explicit stability. The solutions are \Rtwo{computed} with $N=400$ cells. We observe that the numerical investigation of the diffusion and dispersion errors agrees with the properties of the numerical solutions we show. In fact, in Figure~\ref{fig:diss:disp:lintra:p1} we notice that $\gamma=0.5$ provides a very dispersive scheme, whereas $\gamma=1-\sqrt{2}/2$ slightly dampens the extrema of the solution. We also observe that the two third-order schemes, see Figure~\ref{fig:diss:disp:lintra:p2:cou3}, provide almost equivalent approximations in terms of dissipation. However, since they are less dispersive than second-order schemes, they allow to estimate better the speed of propagation and to cure spurious oscillations around discontinuities. 
In conclusion, these results show that the main gain in considering $p=2$ lies in its much lower dispersion compared to $p=1$.

\begin{figure}[t!]
	\centering
	\begin{subfigure}[b]{0.48\textwidth}
		\centering
		\includegraphics[width=\textwidth]{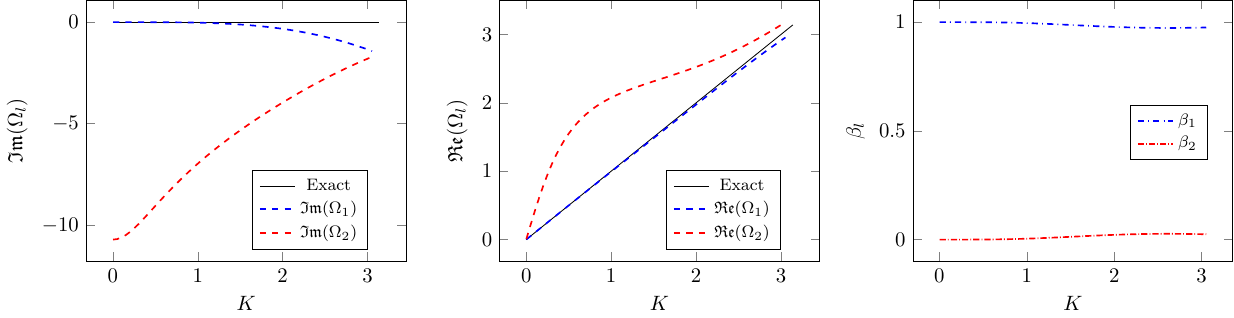}
		\caption{Second-order scheme with $r=1.5$.\label{fig:ev_energy:p1:cou05}}
	\end{subfigure}
	\begin{subfigure}[b]{0.48\textwidth}
		\centering
		\includegraphics[width=\textwidth]{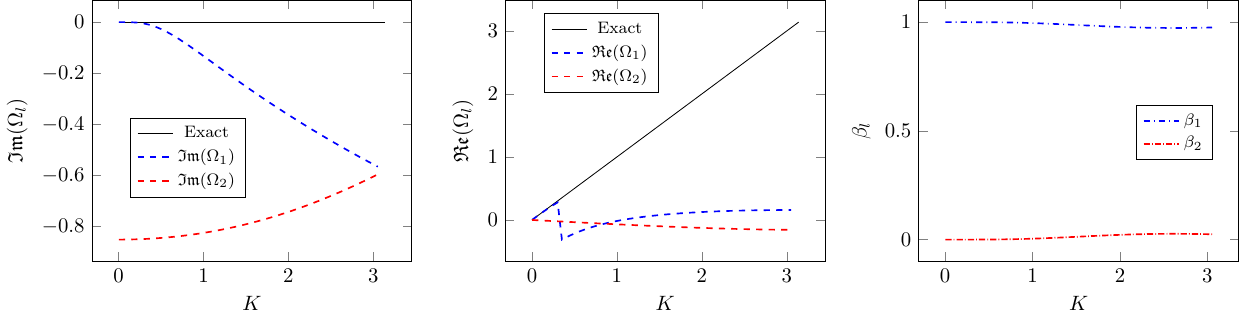}
		\caption{Second-order scheme with $r=15$.\label{fig:ev_energy:p1:cou5}}
	\end{subfigure}
	\begin{subfigure}[b]{0.48\textwidth}
		\centering
		\includegraphics[width=\textwidth]{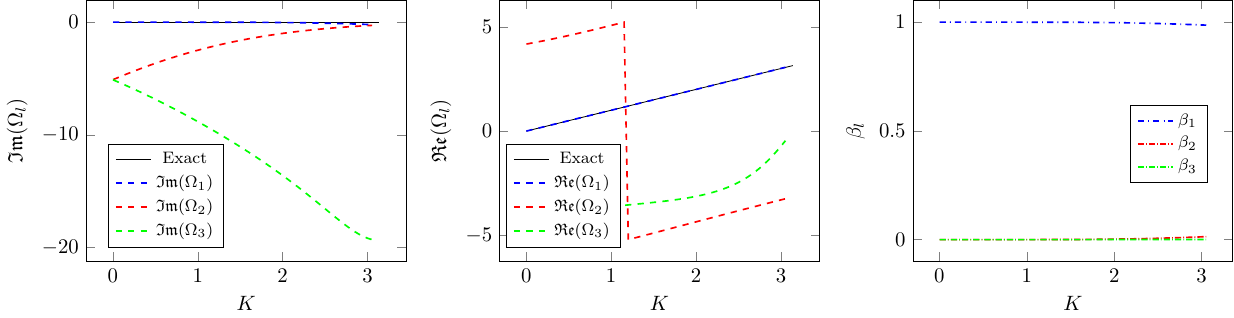}
		\caption{Third-order scheme with $r=1.5$.\label{fig:ev_energy:p2:cou03}}
	\end{subfigure}
	\begin{subfigure}[b]{0.48\textwidth}
		\centering
		\includegraphics[width=\textwidth]{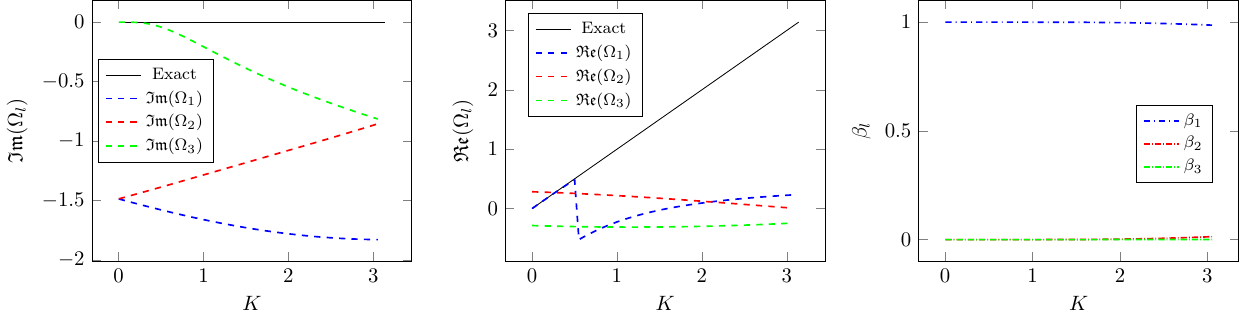}
		\caption{Third-order scheme with $r=15$.\label{fig:ev_energy:p2:cou3}}
	\end{subfigure}
	\caption{Eigenvalues and related energy for $p=1$ (top row) and $p=2$ (bottom row) with DIRK parameters $\gamma=1-\sqrt{2}/2$ and $\gamma=0.4358665215$, respectively.\label{fig:ev_energy}}
\end{figure}

In Figure~\ref{fig:ev_energy}, we illustrate the real and imaginary parts~\eqref{eq:full-discr-omega} of the eigenvalues of the fully-discrete DIRK-DG approximation matrix along with the relative energy. The relative energy has been previously defined, see~\eqref{eq:energy}, and for second-order schemes, it has been analytically proven that as $h\to0$, $\beta_1\to1$ and $\beta_2\to0$. This implies that one eigenvalue becomes dominant. Figure~\ref{fig:ev_energy:p1:cou05} and Figure~\ref{fig:ev_energy:p1:cou5} verify this numerically for second-order schemes, confirming the theoretical results. For third-order schemes, although we were not able to provide a formal proof, Figure~\ref{fig:ev_energy:p2:cou03} and Figure~\ref{fig:ev_energy:p2:cou3} suggest that a similar behavior extends. Specifically, as $h\to0$, the principal eigenvalue remains accurate for larger wave numbers, and the relative energies behave in a corresponding manner.

\section{Space limiting in DIRK-DG schemes} \label{sec:limiting}

In DG approximations, the solution is represented by piecewise polynomial functions defined on each element of the computational domain. At element interfaces, the solution may suffer from spurious oscillations, especially near discontinuities or strong gradients. These oscillations can lead to inaccurate
solutions and even numerical instability. To circumvent the emergence of these undesired oscillations, primarily induced by the high-order spatial approximation, DG discretization incorporates localized limiters. These limiters act upon the high-order coefficients subsequent to each stage of the Runge-Kutta scheme, preventing oscillations while still preserving the accuracy of the solution. The fundamental principle of limiting entails a two-step process: initially, an indicator discerns problematic cells, where local polynomials necessitate correction. Subsequently, within these identified cells, polynomials undergo replacement with modified (limited) counterparts:
\begin{equation}\label{eq:uh:limited}
	\mathbf{u}_{h,\text{lim}}(x,t) = \mathbf{u}^{(0)}_j(t) + \sum_{l=1}^p \mathbf{u}^{(l)}_{j,\text{lim}}(t)\, v^{(l)}_j(x) \quad \mbox{for } x\in I_j,
\end{equation}
where $\mathbf{u}^{(l)}_{j,\text{lim}}(t)$ are the modified moments due to the space limiting procedure.
This ensures that local mean values remain unaltered, thereby upholding conservation principles. The limited polynomial~\eqref{eq:uh:limited} is used to compute the point-values of the solution within the approximation of the space derivative, so that~\eqref{eq:DG_DIRK} becomes
\begin{equation} \label{eq:DG_DIRK:limited}
	\begin{aligned}
		\mathbf{u}^{(l),n+1}_j &= \mathbf{u}^{(l),n}_{j,\text{lim}} - (2l+1)\frac{\Delta t}{h} \sum_{i=1}^s b_i K_{j,\text{lim}}^{(i)} \\
		K_{j,\text{lim}}^{(i)} &= - Q\left(\mathbf{u}_{h,\text{lim}}(\cdot,t^{(i)});l\right) + \left[ \Fnum(\mathbf{u}^{-,(i)}_{j+\frac12,\text{lim}},\mathbf{u}^{+,(i)}_{j+\frac12,\text{lim}}) - (-1)^l \Fnum(\mathbf{u}^{-,(i)}_{j-\frac12,\text{lim}},\mathbf{u}^{+,(i)}_{j-\frac12,\text{lim}}) \right],
	\end{aligned}
\end{equation}
where $\mathbf{u}^{\pm,(i)}_{j+\frac12,\text{lim}}$ are the left and right limits of the limited discontinuous solution $\mathbf{u}_{h,\text{lim}}$ at the cell interface $x_{j+\frac12}$ and time $t^{(i)} = (n + c_i) \Delta t$, and
\begin{equation} \label{eq:integral:discr:limited}
	Q\left(\mathbf{u}_{h,\text{lim}}(\cdot,t);l\right) \approx \disp\int_{-1}^1 \mathbf{f}(\mathbf{u}_{h,\text{lim}}(x_j + \frac{h}{2}y,t))\disp\frac{\mathrm{d} \zeta_l(y)}{\mathrm{d}y} \mathrm{d}y.
\end{equation}

Contemporary RKDG methods offer a diverse array of limiters~\cite{Zhu2013}. As space limiters are typically applied component-wise, we will delineate the two limiters utilized in this paper with respect to a component $u$ within the vector solution $\mathbf{u}$.

\begin{description}	
	\item[Moment limiter of~\cite{BiswasDevineFlaherty94} (BDF limiter)] The moments are modified as follows:
	\begin{equation} \label{eq:limiterBDF}
		u_{j,\text{lim}}^{(l)} = \frac{1}{2l-1} \text{m}\left((2l-1)u_j^{(l)},u_{j+1}^{(l-1)}-u_j^{(l-1)},u_j^{(l-1)}-u_{j-1}^{(l-1)}\right), \quad 1 \leq l \leq p,
	\end{equation}
	where $\text{m}$ is the minmod function given by
	$$
	\text{m}(a_1,a_2,\dots,a_n) = \begin{cases}
		s \cdot \min_{1\leq k \leq n} |a_k|, & \text{if } \text{sign}(a_1)=\dots=\text{sign}(a_n)=s,\\
		0, & \text{otherwise}.
	\end{cases}
	$$
	This limiter is applied adaptively. Initially, the highest-order moment $u_j^{(p)}$ is limited. Subsequently, the limiter is applied to successively lower-order moments whenever the limiting affects the next higher-order moment.
	
	\item[Moment limiter of~\cite{Rider2001} (MP limiter)] The moments are modified as follows:
	\begin{equation} \label{eq:limiterMP}
		u_{j,\text{lim}}^{(l)} = \phi_j u_j^{(l)}, \quad 1 \leq l \leq p,
	\end{equation}
	where
	\begin{equation} \label{eq:phiMP}
	\phi_j = \min\left( 1 , \frac{\Delta u_j^{\max}}{\Delta_{\max} u_j} , \frac{\Delta u_j^{\min}}{\Delta_{\min} u_j} \right).
	\end{equation}
	Here, $\Delta u_j^{\max} = u_j^{\max} - u_j^{(0)}$, where $u_j^{\max}$ is the maximum of the local data, namely $u_j^{\max} = \max(u_{j-1}^{(0)},u_j^{(0)},u_{j+1}^{(0)})$. Further, $\Delta_{\max} u_j = \max_{x\in I_j} u_h(x,t) - u_j^{(0)}$ is the maximum values within the cell $j$ without space limiting minus the cell-average. The MP limiter computes $\max_{x\in I_j} u_h(x,t)$ on a linear reconstruction of the solution: in this case the maximum value occurs at a cell edge, i.e.~$\max_{x\in I_j} u_h(x,t) = \max(u_h(x_{j+1/2},t),u_h(x_{j-1/2},t))$. Indeed, the MP limiter computes $\phi_j$ for the linear reconstruction and then it is used to limit also the higher-order moments. The terms $\Delta u_j^{\min}$ and $\Delta_{\min} u_j$ are computed analogously.
\end{description}

Compared to~\eqref{eq:DIRK:stage}, the space limiting makes system
\begin{equation} \label{eq:DIRK:stage:limited}
	\mathbf{G}_j^{(l)} := \mathbf{u}_j^{(l),(i)} 
	+ (2l+1)\frac{a_{ii}\Delta t}{h} K_{j,\text{lim}}^{(i)} 
	- \mathbf{u}_{j,\text{lim}}^{(l),n} 
	+ (2l+1)\frac{\Delta t}{h} \sum_{\kappa=1}^{i-1} a_{ki} K_{j,\text{lim}}^{(\kappa)} = \mathbf{0}, \quad j=1,\dots,N, \ l=0,\dots,p,
\end{equation}
nonlinear even on linear equations, because of the nonlinear structure of the limiter, which enters in the approximation of the solution at the cell interfaces and at the quadrature nodes.

\subsection{Freezing the nonlinearity with a low-order approximation of the space limiters} \label{sec:freezing}

We propose here an approach to circumvent the nonlinearity due to the space limiters, and keep the nonlinearity of the flux function $\mathbf{f}$ only.

The main idea is to exploit a \emph{predictor} $\mathbf{u}^{\star}_j$ of the solution to pre-compute and freeze the nonlinear space limiters in~\eqref{eq:limiterMP} so that the limited DG moments are
\begin{equation} \label{eq:limiterMP:pred}
	\begin{aligned}
	u_{j,\text{lim}^\star}^{(l)} &= \phi_j^\star u_j^{(l)}, \quad 1 \leq l \leq p,\\
	\phi_j^\star &= \min\left( 1 , \frac{\Delta u_j^{\max}}{\Delta_{\max} u_j} , \frac{\Delta u_j^{\min}}{\Delta_{\min} u_j} \right).
	\end{aligned}
\end{equation}
Observe that $\phi_j^\star$ is as in~\eqref{eq:phiMP}, but it is now computed on the predictor solution. Namely, one has that $\Delta u_j^{\max} = u_j^{\max,\star} - u_j^{\star}$, with $u_j^{\max,\star} = \max(u_{j-1}^{\star},u_j^{\star},u_{j+1}^{\star})$, and $\Delta_{\max} u_j = \max_{x\in I_j} u_h^{\star}(x,t) - u_j^{\star}$, with $\max_{x\in I_j} u_h^{\star}(x,t) = \max(u_h^{\star}(x_{j+1/2},t),u_h^{\star}(x_{j-1/2},t))$ and
$$
	u_h^{\star}(x_{j+1/2},t) = u_j^{\star}(t) + \frac{u_{j+1}^{\star}(t)-u_j^{\star}(t)}{2}, \quad
	u_h^{\star}(x_{j-1/2},t) = u_j^{\star}(t) - \frac{u_j^{\star}(t)-u_{j-1}^{\star}(t)}{2}.
$$
The terms $\Delta u_j^{\min}$ and $\Delta_{\min} u_j$ are also computed on the predictor solution $\mathbf{u}^\star$. \Rtwo{Specifically, $\Delta u_j^{\min}=u_j^{\star}-\min(u_{j-1}^{\star},u_j^{\star},u_{j+1}^{\star})$ and $\Delta_{\min} u_j = u_j^{\star} - \min(u_h^{\star}(x_{j+1/2},t),u_h^{\star}(x_{j-1/2},t))$}. In this way, the complete scheme would be linear with respect to the space limiters, and nonlinear only through the flux function. With this definition, $\phi_j^\star \in \{0,1\}$ and therefore it allows us to detect whether a cell is troubling or not.

The idea of freezing the nonlinearity introduced by the high-order space approximation relies on the works~\cite{Puppo2023,Puppo2024} where it was devised in the finite-volume framework. Here, we explore the same approach for high-order DG approximations. First, we describe the choice of the predictor scheme and then we formulate the final \emph{semi-linear} high-order scheme.

\subsubsection{The space-time first-order implicit predictor} \label{sec:limiting:predictor}

In order to compute the predictor of the DIRK stages, we consider a $\mathbb{P}_0$ DG approximation and we integrate system~\eqref{eq:DG_odesystem}, with $p=0$, using an implicit first-order scheme. This leads us to a composite backward Euler method with Butcher tableau
\begin{equation}
	\label{eq:tableau:be}
	\begin{array}{c|cccc}
		c_1 & c_1 & 0 & \dots & 0 \\[1.5ex]
		c_2 & c_1 & c_2-c_1 & \dots & 0 \\[1.5ex]
		\vdots & \vdots & \vdots & \ddots & \vdots \\[1.5ex]
		c_s & c_1 & c_2-c_1 & \dots & c_s-c_{s-1} \\[1ex]
		\hline
		&&&&\\[-1.8ex]
		& c_1 & c_2-c_1 & \dots & c_s-c_{s-1}
	\end{array}
\end{equation}
which provides the $s$ first-order approximations $\mathbf{u}_j^{\star,(i)}$, $j=1,\dots,N$, $i=1,\dots,s$, at any time $t^{(i)}=(n+c_{i}) \Delta t \in [n\Delta t,(n+1)\Delta t]$, where $\Delta t$ is the time-step of the high-order DIRK scheme, as follows
\begin{subequations} \label{eq:predictor}
	\begin{align}
		\mathbf{u}_j^{\star,(i)} &= \mathbf{u}_j^{\star,(i-1)} - (c_i-c_{i-1}) \frac{\Delta t}{h} K_j^{\star,(i)}, \label{eq:predictor:stage} \\
		K_j^{\star,(i)} &= \Fnum\left(\mathbf{u}_{j}^{\star,(i)},\mathbf{u}_{j+1}^{\star,(i)}\right) - \Fnum\left(\mathbf{u}_{j-1}^{\star,(i)},\mathbf{u}_{j}^{\star,(i)}\right). \label{eq:predictor:flux}
	\end{align}
\end{subequations}
We use the convention $c_{0}=0$  and, in addition, $\mathbf{u}_j^{\star,(0)}:=\mathbf{u}_j^{(0)}$, where $\mathbf{u}_j^{(0)}$ is the approximation of the zero-th DG moment with the high-order scheme. Instead, the values $\mathbf{u}_{0,N+1}^{\star,(0)}$ are assigned by boundary conditions.

Notice that the numerical flux functions $\Fnum$ in~\eqref{eq:predictor:flux} is computed on piecewise constant, unlimited, approximations of the cell boundary values since first-order schemes do not require space-limiting, because they are unconditionally Total Variation Diminishing. Therefore, the first-order scheme~\eqref{eq:predictor} is characterized by a single nonlinearity, that is the one induced by the flux function $\mathbf{f}$. This means that the computation of~\eqref{eq:predictor:stage} requires the solution of $s$ nonlinear systems of dimension $m N$, which are fully linear with respect to the unknowns $\mathbf{U}^{\star,(i)}=\{\mathbf{u}^{\star,(i)}_j\}_{j=1}^N$ if the flux function $\mathbf{f}$ is linear:
\begin{equation} \label{eq:predictor:system}
	\mathbf{G}(\mathbf{U}^{\star,(i)}) := 
	\mathbf{U}^{\star,(i)} 
	+ (c_i-c_{i-1}) \frac{\Delta t}{h} \mathbf{K}^{\star,(i)} 
	- \mathbf{U}^{\star,(i-1)} = \mathbf{0}, \quad i=1,\dots,s,
\end{equation}
where $\mathbf{K}^{\star,(i)}=\{K^{\star,(i)}_j\}_{j=1}^N$.
The solution of~\eqref{eq:predictor:system} requires the use of a nonlinear solver $s$ times within a single time-step. In~\cite{Puppo2023,Puppo2024}, based on finite volume space discretization, the authors rely on classical Newton-Raphson method:
\begin{equation} \label{eq:newton}
	\frac{\partial \mathbf{G}(\mathbf{U}_{(k)}^{\star,(i)})}{\partial \mathbf{U}^{\star,(i)}} \Delta \mathbf{u} = - \mathbf{G}(\mathbf{U}_{(k)}^{\star,(i)}), \quad \mathbf{U}_{(0)}^{\star,(i)} = \mathbf{U}^{\star,(i-1)}, \quad i=1,\dots,s
\end{equation}
where $\mathbf{U}^{\star,(i)}_{(k)}$ it the Newton's iteration, $k\geq 0$ is the iteration counter, and $\Delta \mathbf{u} = \mathbf{U}^{\star,(i)}_{(k+1)}-\mathbf{U}^{\star,(i)}_{(k)}$.
In the DG frameworks, however, the size of the unknowns is much larger and, thus, the formation of the Jacobian matrix leads to a very large cost. Therefore, to limit the cost due to the formation of the Jacobian, we employ a Jacobian-free Newton-Krylov method~\cite{Knoll2004} whereby the Jacobian $J^{\star,(i)}\in\mathbb{R}^{mN\times mN}$ is approximated by finite difference
$$
	J^{\star,(i)} \approx \frac{\mathbf{G}(\mathbf{U}^{\star,(i)}+\epsilon\Delta \mathbf{u})-\mathbf{G}(\mathbf{U}^{\star,(i)})}{\epsilon \Delta \mathbf{u}},
$$
where $\epsilon$ is a small quantity, \Rtwo{of order $\sim 10^{-5}$}. Direct substitution into classical Newton's iteration~\eqref{eq:newton} leads to the residual problem
$$
	\frac{\mathbf{G}(\mathbf{U}^{\star,(i)}+\epsilon\Delta \mathbf{u})-\mathbf{G}(\mathbf{U}^{\star,(i)})}{\epsilon} = -\mathbf{G}(\mathbf{U}^{\star,(i)})
$$
that is solved iteratively by the GMRES method.

To address potential spurious oscillations at high Courant numbers, the cells identified as problematic by the MP limiter, computed on the predictor, are expanded to include adjacent cells that are also marked as troubling. This extension is contingent upon the Courant number, ensuring that the limiting process accounts for the increased numerical flux influences at larger time-steps. The dependency on the Courant number and its impact on the extension of troubling cells are investigated in the numerical results section. The extension will be denoted by the parameter $\delta$.

\begin{remark}
	We list below the Butcher tableaux of the predictor schemes corresponding to the DIRK methods given in Table~\ref{tab:labels}.
	\begin{table}[h!]
		\begin{minipage}{.25\linewidth}
			\centering
			\footnotesize
			\[\arraycolsep=5.0pt\def\arraystretch{1.4}
			\begin{array}{c|cc}
				\gamma & \gamma & 0 \\[1.5ex]
				1-\gamma & \gamma & 1-2\gamma \\[1.5ex]
				\hline
				& \gamma & 1-2\gamma 
			\end{array}\]
			\caption{Predictor scheme of the $2$-stage DIRK method in Table~\ref{tab:dirk22}.\label{tab:predictor:of:dirk22}}
		\end{minipage}
		\hfill
		\begin{minipage}{0.35\linewidth}
			\centering
			\footnotesize
			\[\arraycolsep=5.0pt\def\arraystretch{1.5}
			\begin{array}{c|ccc}
				\gamma & \gamma & 0 & 0 \\[1.4ex]
				1+\tfrac{\gamma}{2} & \gamma & 1-\tfrac{\gamma}{2} & 0 \\[1.5ex]
				1 & \gamma & 1-\tfrac{\gamma}{2} & -\tfrac{\gamma}{2} \\[1.5ex]
				\hline
				& \gamma & 1-\tfrac{\gamma}{2} & -\tfrac{\gamma}{2}
			\end{array}\]
			\caption{Predictor scheme of the $3$-stage DIRK method in Table~\ref{tab:dirk33}.\label{tab:predictor:of:dirk33}}
		\end{minipage}
		\hfill
		\begin{minipage}{0.25\linewidth}
			\centering
			\footnotesize
			\[\arraycolsep=5.0pt\def\arraystretch{1.5}
			\begin{array}{c|cccc}
				\tfrac12 & \tfrac12 & 0 & 0 & 0 \\[1.4ex]
				\tfrac23 & \tfrac12 & \tfrac16 & 0 & 0 \\[1.5ex]
				\tfrac12 & \tfrac12 & \tfrac16 & -\tfrac16 & 0 \\[1.5ex]
				1 & \tfrac12 & \tfrac16 & -\tfrac16 & \tfrac12 \\[1.5ex]
				\hline
				& \tfrac12 & \tfrac16 & -\tfrac16 & \tfrac12
			\end{array}\]
			\caption{Predictor scheme of the $4$-stage third-order DIRK method in Table~\ref{tab:dirk43}.\label{tab:predictor:of:dirk43}}
		\end{minipage}
	\end{table}
\end{remark}

\subsubsection{The space-time high-order \emph{semi-linearized} implicit scheme}

Assume that the first-order predictor $\mathbf{U}^{\star,(i)}$ is available at time $t^{(i)} \in [n\Delta t,(n+1)\Delta t]$. Then, we use it to estimate the nonlinear MP limiter as explained at the beginning of Section~\ref{sec:freezing}, see~\eqref{eq:limiterMP:pred}, which gives us the limited (on the predictor) discontinuous solution
\begin{equation}\label{eq:uh:limited:pred}
	\mathbf{u}_{h,\text{lim}^\star}(x,t) = \mathbf{u}^{(0)}_j(t) + \sum_{l=1}^p \mathbf{u}^{(l)}_{j,\text{lim}^\star}(t)\, v^{(l)}_j(x) \quad \mbox{for } x\in I_j.
\end{equation}
Therefore, the high-order solution at time $t^{(i)}$ can be computed as in~\eqref{eq:DG_DIRK:limited} with
\begin{equation} \label{eq:spaceder:limited:pred}
	K_{j,\text{lim}^\star}^{(i)} = - Q\left(\mathbf{u}_{h,\text{lim}^\star}(\cdot,t^{(i)});l\right) + \left[ \Fnum(\mathbf{u}^{-,(i)}_{j+\frac12,\text{lim}^\star},\mathbf{u}^{+,(i)}_{j+\frac12,\text{lim}^\star}) - (-1)^l \Fnum(\mathbf{u}^{-,(i)}_{j-\frac12,\text{lim}^\star},\mathbf{u}^{+,(i)}_{j-\frac12,\text{lim}^\star}) \right],
\end{equation}
which is the space approximation at the $i$-th DIRK stage. We use the subscript $\text{lim}^\star$ to stress that this is computed by freezing the limiters on the predictor solution. This term is now fully linear with respect to the unknowns at time $t^{(i)}=t_n+c_i\Delta t$, unless the non-linearity of the flux function $\mathbf{f}$. Therefore, we approximate the nonlinear system~\eqref{eq:DIRK:stage:limited} with
\begin{equation} \label{eq:DIRK:stage:limited:pred}
	\mathbf{G}_j^{(l)} := \mathbf{u}_j^{(l),(i)} 
	+ (2l+1)\frac{a_{ii}\Delta t}{h} K_{j,\text{lim}^\star}^{(i)} 
	- \mathbf{u}_{j,\text{lim}}^{(l),n} 
	+ (2l+1)\frac{\Delta t}{h} \sum_{\kappa=1}^{i-1} a_{ki} K_{j,\text{lim}}^{(\kappa)} = \mathbf{0}, \quad j=1,\dots,N, \ l=0,\dots,p,
\end{equation}
that is nonlinear only through the nonlinear physical flux. 

Once the stage values $\mathbf{U}^{(l),(i)}=\{\mathbf{u}_j^{(l),(i)}\}_{j=1}^N$, $l=0,\dots,p$, are computed, they still undergo the limiting process using the BDF limiter~\eqref{eq:limiterBDF}. In fact, while the MP limiter is suited for the predictor, due to the structure~\eqref{eq:limiterMP}, the BDF one limits adaptively the higher-order moments allowing for a better accuracy. The resulting limited stage values $\mathbf{U}_{\text{lim}}^{(l),(i)}$ are then used to compute $\mathbf{K}_{\text{lim}}^{(i)}=\{K_{j,\text{lim}}^{(i)}\}_{j=1}^N$ which enters in the computation of the successive stage values $\mathbf{U}_{\text{lim}}^{(l),(\kappa)}$, $\kappa>i$.

System~\eqref{eq:DIRK:stage:limited:pred} is also numerically tackled using the Jacobian-free Newton-Krylov method described previously.

Finally, the high-order solution $\mathbf{U}^{(l),n+1}=\{\mathbf{u}_j^{(l),n+1}\}_{j=1}^N$, $l=0,\dots,p$, at time level $(n+1)\Delta t$ is obtained as in~\eqref{eq:DG_DIRK:limited}. 
Overall, the computation of the solution with the DIRK-DG scheme is sketched in Algorithm~\ref{alg:dirkdg}.

\begin{algorithm}
	\caption{$s$-stage DIRK-DG of order $p+1$}\label{alg:dirkdg}
	\begin{algorithmic}[1]
		\INPUT High-order limited solution $\mathbf{u}_{j,\text{lim}}^{(l),n}$ at time $t^n$, for $j=1,\dots,N$, $l=0,\dots,p$
		\FOR{$i = 1$ to $s$}
		\STATE Solve~\eqref{eq:predictor:system} to get the predictor solution $\mathbf{u}_j^{\star,(i)}$, $j=1\dots,N$, at time $t^{(i)}=t^n+c_i\Delta t$
		\STATE Compute the MP limiters $\phi_j^*$, $j=1,\dots,N$, as in~\eqref{eq:limiterMP:pred}
		\IF{$\phi_j^* = 0$ for some $j=1,\dots,N$}
		\STATE impose $\phi_j^* = 0$ also for the $2\delta$ adjacent cells $\{I_k\}_{k=j-\delta,k\neq j}^{j+\delta}$
		\ENDIF
		\STATE Solve~\eqref{eq:DIRK:stage:limited:pred} to get the high-order solution $\mathbf{u}_j^{(l),(i)}$, $j=1\dots,N$, $l=0,\dots,p$, at time $t^{(i)}$
		\STATE Limit the high-order stage $\mathbf{u}_j^{(l),(i)}$ with the BDF limiter~\eqref{eq:limiterBDF}, to get the limited solution $\mathbf{u}_{j,\text{lim}}^{(l),(i)}$
		\ENDFOR
		\STATE Compute the high-order limited solution $\mathbf{u}_{j,\text{lim}}^{(l),n+1}$ at time $t^{n+1}$, for $j=1,\dots,N$, $l=0,\dots,p$, as in~\eqref{eq:DG_DIRK:limited}
	\end{algorithmic}
\end{algorithm}

\section{Numerical simulations} \label{sec:numerical:simulations}

For simplicity and not to distract from the main focus of the paper, which is on the implicit time-integration of high-order DG schemes, all the numerical simulations are performed with the Rusanov (local Lax-Friedrichs) numerical flux
\begin{equation}\label{eq:LFflux}
	\Fnum(\mathbf{u}^-,\mathbf{u}^+) = \frac{ \mathbf{f}(\mathbf{u}^-)+\mathbf{f}(\mathbf{u}^+) }{2} - \frac{\alpha(\mathbf{u}^+,\mathbf{u}^-)}{2} (\mathbf{u}^+-\mathbf{u}^-),
\end{equation}
where the viscosity parameter $\alpha$ will be discussed for each numerical example. In explicit schemes one has to choose $\alpha = \max\{\|\mathbf{f}'(\mathbf{v})\|,\|\mathbf{f}'(\mathbf{w})\|\}$, where we recall that $\|\mathbf{f}'(\cdot)\|$ denotes the spectral radius of the Jacobian of the flux function $\mathbf{f}$. In implicit schemes this choice would deserve more attention, especially if the fluid speed is much lower than the sound speed and if one is not interested in the acoustic waves. In this case one can use more appropriate speed estimates in the approximate Riemann solver as explained, e.g., in~\cite{Degond2011}.

All numerical simulations are performed with the DIRK-DG schemes listed in Table~\ref{tab:labels} and described by Algorithm~\ref{alg:dirkdg}. In the case of hyperbolic systems, they will be compared with explicit RKDG schemes and the implicit Quinpi finite-volume scheme of~\cite{Puppo2023,Puppo2024}. We recall that for the solution of the nonlinear systems in the predictor and in the high-order steps, we use the Jacobian-Free Newton-Krylov method with a simple GMRES linear solver without preconditioner~\cite{Knoll2004}. For a detailed discussion on properties of iterative solvers in the context of implicit discretizations of hyperbolic equations we refer to~\cite{Birken2022}.

Finally, we recall the definition of the ratio $r$ between the time-step used and the time-step required for explicit stability, as defined in equation~\eqref{eq:def:Courant}. This parameter will be extensively used in the following sections to describe the regimes under which the schemes operate.


\subsection{Scalar conservation laws}

The following scalar tests are intended to serve as a proof of concept for the implicit schemes, therefore we are particularly interested in regimes where $r>1$. See~\eqref{eq:def:Courant} for the definition of $r$. Later, when considering systems of conservation laws, the value of $r$ will be precisely determined based on the accuracy requirements for material waves.

For scalar conservation laws, the viscosity parameter of the Rusanov flux is obviously chosen as $\alpha(u^+,u^-) = \max \{ |f'(u^+)|,|f'(u^-)| \}$.

\subsubsection{Convergence test}

We test the numerical convergence rate of DIRK-DG schemes on the linear scalar conservation law~\eqref{eq:adv_lin} with speed $a=1$ on $[-1,1]$ with periodic boundary conditions, up to the final time $t = 2$. The initial condition
\begin{equation} \label{eq:convtestIC}
	u_0(x) = \sin(\pi x - \frac{1}{\pi} \sin(\pi x))
\end{equation}
is taken from~\cite{Arandiga2011}. The numerical errors in $L^1$ norm and convergence rates are showed in Table~\ref{tab:rates:p1:DIRK:025}, Table~\ref{tab:rates:p1:DIRK:Lstable} and Table~\ref{tab:rates:p1:DIRK:05} for second-order DIRK-DG scheme, and in Table~\ref{tab:rates:p2:DIRK:Lstable} and Table~\ref{tab:rates:p2:DIRK43} for third-order DIRK-DG schemes. All the errors and experimental orders are computed for different Courant numbers.

\begin{table}[th!]
	\caption{Experimental order of convergence of the second-order scheme $\DIRKdue$, computed on the initial condition~\eqref{eq:convtestIC} and the DIRK scheme in Table~\ref{tab:dirk22} with $\gamma=0.25$.\label{tab:rates:p1:DIRK:025}}
	\centering
	\vspace{0.25cm}
	\pgfplotstabletypeset[
	font=\small,
	col sep=comma,
	sci zerofill,
	empty cells with={--},
	every head row/.style={before row={\toprule
			&
			&\multicolumn{2}{c}{$r=1$}
			&\multicolumn{2}{c}{$r=3$}
			&\multicolumn{2}{c}{$r=9$}
			&\multicolumn{2}{c}{$r=15$}
			&\multicolumn{2}{c}{$r=30$}
			\\
		},
		after row=\midrule
	},
	every last row/.style={after row=\bottomrule},
	create on use/rate1/.style={create col/dyadic refinement rate={1}},
	create on use/rate2/.style={create col/dyadic refinement rate={2}},
	create on use/rate3/.style={create col/dyadic refinement rate={3}},
	create on use/rate4/.style={create col/dyadic refinement rate={4}},
	create on use/rate5/.style={create col/dyadic refinement rate={5}},
	columns/0/.style={column name={$N$}},
	columns/6/.style={column name={$h$}},
	columns/1/.style={column name={$L^1$ error},sci e},
	columns/rate1/.style={fixed zerofill,column name={rate}},
	columns/2/.style={column name={$L^1$ error},sci e},
	columns/rate2/.style={fixed zerofill,column name={rate}},
	columns/3/.style={column name={$L^1$ error},sci e},
	columns/rate3/.style={fixed zerofill,column name={rate}},
	columns/4/.style={column name={$L^1$ error},sci e},
	columns/rate4/.style={fixed zerofill,column name={rate}},
	columns/5/.style={column name={$L^1$ error},sci e},
	columns/rate5/.style={fixed zerofill,column name={rate}},
	columns={6,0,1,rate1,2,rate2,3,rate3,4,rate4,5,rate5},
	]
	{err_p1_DIRK_025.err}
\end{table}

\begin{table}[th!]
	\caption{Experimental order of convergence of the second-order scheme $\DIRKlstab$, computed on the initial condition~\eqref{eq:convtestIC} and the DIRK scheme in Table~\ref{tab:dirk22} with $\gamma=1-\sqrt{2}/2$.\label{tab:rates:p1:DIRK:Lstable}}
	\centering
	\vspace{0.25cm}
	\pgfplotstabletypeset[
	font=\small,
	col sep=comma,
	sci zerofill,
	empty cells with={--},
	every head row/.style={before row={\toprule
			&
			&\multicolumn{2}{c}{$r=1$}
			&\multicolumn{2}{c}{$r=3$}
			&\multicolumn{2}{c}{$r=9$}
			&\multicolumn{2}{c}{$r=15$}
			&\multicolumn{2}{c}{$r=30$}
			\\
		},
		after row=\midrule
	},
	every last row/.style={after row=\bottomrule},
	create on use/rate1/.style={create col/dyadic refinement rate={1}},
	create on use/rate2/.style={create col/dyadic refinement rate={2}},
	create on use/rate3/.style={create col/dyadic refinement rate={3}},
	create on use/rate4/.style={create col/dyadic refinement rate={4}},
	create on use/rate5/.style={create col/dyadic refinement rate={5}},
	columns/0/.style={column name={$N$}},
	columns/6/.style={column name={$h$}},
	columns/1/.style={column name={$L^1$ error},sci e},
	columns/rate1/.style={fixed zerofill,column name={rate}},
	columns/2/.style={column name={$L^1$ error},sci e},
	columns/rate2/.style={fixed zerofill,column name={rate}},
	columns/3/.style={column name={$L^1$ error},sci e},
	columns/rate3/.style={fixed zerofill,column name={rate}},
	columns/4/.style={column name={$L^1$ error},sci e},
	columns/rate4/.style={fixed zerofill,column name={rate}},
	columns/5/.style={column name={$L^1$ error},sci e},
	columns/rate5/.style={fixed zerofill,column name={rate}},
	columns={6,0,1,rate1,2,rate2,3,rate3,4,rate4,5,rate5},
	]
	{err_p1_DIRK_Lstable.err}
\end{table}

\begin{table}[th!]
	\caption{Experimental order of convergence of the second-order scheme $\DIRKlam{0.5}$, computed on the initial condition~\eqref{eq:convtestIC} and the DIRK scheme in Table~\ref{tab:dirk22} with $\gamma=0.5$.\label{tab:rates:p1:DIRK:05}}
	\centering
	\vspace{0.25cm}
	\pgfplotstabletypeset[
	font=\small,
	col sep=comma,
	sci zerofill,
	empty cells with={--},
	every head row/.style={before row={\toprule
			&
			&\multicolumn{2}{c}{$r=1$}
			&\multicolumn{2}{c}{$r=3$}
			&\multicolumn{2}{c}{$r=9$}
			&\multicolumn{2}{c}{$r=15$}
			&\multicolumn{2}{c}{$r=30$}
			\\
		},
		after row=\midrule
	},
	every last row/.style={after row=\bottomrule},
	create on use/rate1/.style={create col/dyadic refinement rate={1}},
	create on use/rate2/.style={create col/dyadic refinement rate={2}},
	create on use/rate3/.style={create col/dyadic refinement rate={3}},
	create on use/rate4/.style={create col/dyadic refinement rate={4}},
	create on use/rate5/.style={create col/dyadic refinement rate={5}},
	columns/0/.style={column name={$N$}},
	columns/6/.style={column name={$h$}},
	columns/1/.style={column name={$L^1$ error},sci e},
	columns/rate1/.style={fixed zerofill,column name={rate}},
	columns/2/.style={column name={$L^1$ error},sci e},
	columns/rate2/.style={fixed zerofill,column name={rate}},
	columns/3/.style={column name={$L^1$ error},sci e},
	columns/rate3/.style={fixed zerofill,column name={rate}},
	columns/4/.style={column name={$L^1$ error},sci e},
	columns/rate4/.style={fixed zerofill,column name={rate}},
	columns/5/.style={column name={$L^1$ error},sci e},
	columns/rate5/.style={fixed zerofill,column name={rate}},
	columns={6,0,1,rate1,2,rate2,3,rate3,4,rate4,5,rate5},
	]
	{err_p1_DIRK_05.err}
\end{table}

\begin{table}[th!]
	\caption{Experimental order of convergence of the third-order scheme $\DIRKlstabTre$, computed on the initial condition~\eqref{eq:convtestIC} and the DIRK scheme in Table~\ref{tab:dirk33} with $\gamma=0.435866521508459$.\label{tab:rates:p2:DIRK:Lstable}}
	\centering
	\vspace{0.25cm}
	\pgfplotstabletypeset[
	font=\small,
	col sep=comma,
	sci zerofill,
	empty cells with={--},
	every head row/.style={before row={\toprule
			&
			&\multicolumn{2}{c}{$r=1$}
			&\multicolumn{2}{c}{$r=5$}
			&\multicolumn{2}{c}{$r=15$}
			&\multicolumn{2}{c}{$r=25$}
			&\multicolumn{2}{c}{$r=50$}
			\\
		},
		after row=\midrule
	},
	every last row/.style={after row=\bottomrule},
	create on use/rate1/.style={create col/dyadic refinement rate={1}},
	create on use/rate2/.style={create col/dyadic refinement rate={2}},
	create on use/rate3/.style={create col/dyadic refinement rate={3}},
	create on use/rate4/.style={create col/dyadic refinement rate={4}},
	create on use/rate5/.style={create col/dyadic refinement rate={5}},
	columns/0/.style={column name={$N$}},
	columns/6/.style={column name={$h$}},
	columns/1/.style={column name={$L^1$ error},sci e},
	columns/rate1/.style={fixed zerofill,column name={rate}},
	columns/2/.style={column name={$L^1$ error},sci e},
	columns/rate2/.style={fixed zerofill,column name={rate}},
	columns/3/.style={column name={$L^1$ error},sci e},
	columns/rate3/.style={fixed zerofill,column name={rate}},
	columns/4/.style={column name={$L^1$ error},sci e},
	columns/rate4/.style={fixed zerofill,column name={rate}},
	columns/5/.style={column name={$L^1$ error},sci e},
	columns/rate5/.style={fixed zerofill,column name={rate}},
	columns={6,0,1,rate1,2,rate2,3,rate3,4,rate4,5,rate5},
	]
	{err_p2_DIRK_Lstable.err}
\end{table}

\begin{table}[th!]
	\caption{Experimental order of convergence of the third-order scheme $\DIRKQuattroTre$, computed on the initial condition~\eqref{eq:convtestIC} and the four-stage DIRK scheme in Table~\ref{tab:dirk43}.\label{tab:rates:p2:DIRK43}}
	\centering
	\vspace{0.25cm}
	\pgfplotstabletypeset[
	font=\small,
	col sep=comma,
	sci zerofill,
	empty cells with={--},
	every head row/.style={before row={\toprule
			&
			&\multicolumn{2}{c}{$r=1$}
			&\multicolumn{2}{c}{$r=5$}
			&\multicolumn{2}{c}{$r=15$}
			&\multicolumn{2}{c}{$r=25$}
			&\multicolumn{2}{c}{$r=50$}
			\\
		},
		after row=\midrule
	},
	every last row/.style={after row=\bottomrule},
	create on use/rate1/.style={create col/dyadic refinement rate={1}},
	create on use/rate2/.style={create col/dyadic refinement rate={2}},
	create on use/rate3/.style={create col/dyadic refinement rate={3}},
	create on use/rate4/.style={create col/dyadic refinement rate={4}},
	create on use/rate5/.style={create col/dyadic refinement rate={5}},
	columns/0/.style={column name={$N$}},
	columns/6/.style={column name={$h$}},
	columns/1/.style={column name={$L^1$ error},sci e},
	columns/rate1/.style={fixed zerofill,column name={rate}},
	columns/2/.style={column name={$L^1$ error},sci e},
	columns/rate2/.style={fixed zerofill,column name={rate}},
	columns/3/.style={column name={$L^1$ error},sci e},
	columns/rate3/.style={fixed zerofill,column name={rate}},
	columns/4/.style={column name={$L^1$ error},sci e},
	columns/rate4/.style={fixed zerofill,column name={rate}},
	columns/5/.style={column name={$L^1$ error},sci e},
	columns/rate5/.style={fixed zerofill,column name={rate}},
	columns={6,0,1,rate1,2,rate2,3,rate3,4,rate4,5,rate5},
	]
	{err_p2_DIRK43.err}
\end{table}

As anticipated by the Fourier analysis in Section~\ref{sec:analysis}, the implicit scheme $\DIRKlam{0.5}$ is the least accurate among the tested second-order schemes, due to its higher dispersive error. For this reason it will not be considered anymore in the following. The third-order scheme $\DIRKQuattroTre$ is more accurate than $\DIRKlstabTre$. All the tested schemes achieve the theoretical order of convergence even with relatively coarse grids. This indicates that the presence of the predictor does not introduce unnecessary limiting, as expected for a smooth profile. Therefore, at the same order, the differences in results are solely due to the different time-integration schemes. Finally, we notice that as the ratio between the explicit time-step and the effective time-step increases, the errors increase, as expected.

In Table~\ref{tab:rates:p2:DIRK:Lstable} and Table~\ref{tab:rates:p2:DIRK43}, we computed the errors and rates for $N=640$ and $N=1280$ cells for the case $r=5$. This was done to demonstrate that the experimental order of convergence has a local decrease around $N=320$ cells, but on average it is third-order accurate.

\subsubsection{Linear advection}

We consider again the linear scalar conservation law~\eqref{eq:adv_lin} with speed $a=1$ on the periodic domain $x\in[-1,1]$. Now, we evolve initial discontinuous profiles $u_0(x)$ for one period, i.e.~up to final time $t=2$. Precisely, we consider the following non-smooth initial conditions
	\begin{equation} \label{eq:doublestepIC}
		u_0(x) =
		\begin{cases}
			1, & -0.25 \leq x \leq 0.25,\\
			0, & \text{otherwise.}
		\end{cases}
	\end{equation}
and
	\begin{equation} \label{eq:sindiscontIC}
		u_0(x) = \sin(\pi x) +
		\begin{cases}
			3, & -0.4 \leq x \leq 0.4,\\
			0, & \text{otherwise,}	
		\end{cases}
	\end{equation}
In Section~\ref{sec:analysis}, we numerically studied the properties of DIRK-DG schemes for transporting non-smooth data with minimal dissipation and dispersion effects. However, we observed that the use of spatial limiters is necessary. Therefore, in this section, we investigate the effectiveness of the predictor in detecting spurious oscillations near discontinuities. 

\begin{figure}[t!]
	\centering
	\begin{subfigure}[b]{0.48\textwidth}
		\centering
		\includegraphics[width=\textwidth]{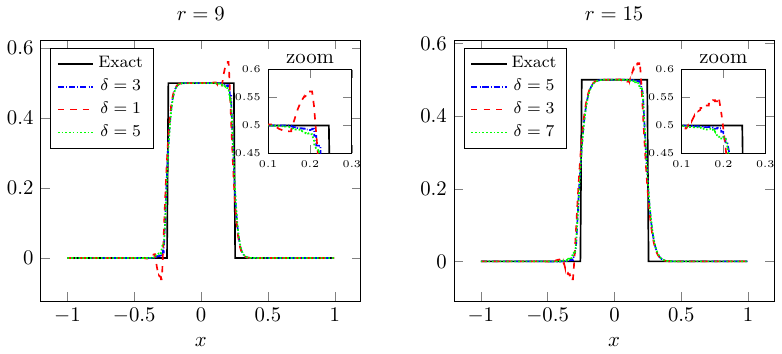}
		\caption{Second-order scheme $\DIRKdue$.\label{fig:li:doubleStep:p1:param:ext}}
	\end{subfigure}
	\hfill
	\begin{subfigure}[b]{0.48\textwidth}
		\centering
		\includegraphics[width=\textwidth]{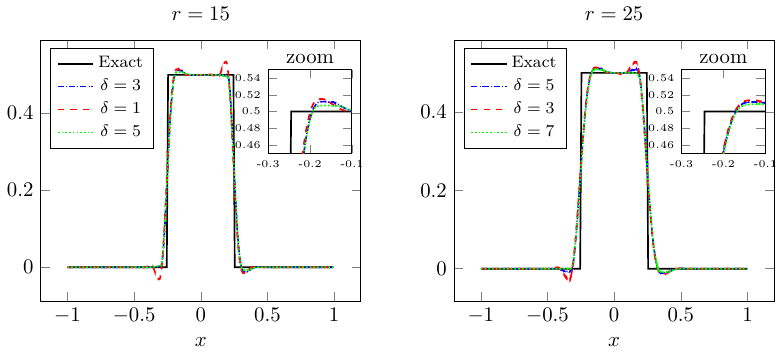}
		\caption{Third-order scheme $\DIRKlstabTre$.\label{fig:li:doubleStep:p2:param:ext}}
	\end{subfigure}
	\caption{Solutions to the linear advection equation with speed $a=1$, final time $t=2$ and initial condition~\eqref{eq:doublestepIC}. The numerical approximations are obtained with $N=400$ cells. The different solutions show the effect of the extension parameter $\delta$ within the space limiting based on the first-order predictor.\label{fig:li:doubleStep:param:ext}}
\end{figure}

First of all, to deal with the problem of oscillations, Figure~\ref{fig:li:doubleStep:param:ext} shows the effect of the \textit{extension parameter} $\delta$ introduced at the end of Section~\ref{sec:limiting:predictor}. We recall that $\delta$ defines the number of cells marked as troubling to the left and right of a cell identified as non smooth by the MP limiter. See also Step 5 in Algorithm~\ref{alg:dirkdg}. We use $r=9,15$ for second-order schemes and $r=15,25$ for third-order schemes. Precisely, in order to test the role of the parameter $\delta$ we use the methods $\DIRKdue$, see Figure~\ref{fig:li:doubleStep:p1:param:ext}, and $\DIRKlstabTre$, see Figure~\ref{fig:li:doubleStep:p2:param:ext}. The optimal value of $\delta$ turns out to be $\delta\approx r/(2p+1)$. When it is too small, the fluctuations are not smoothed, while when it is too large, it has no particular effect on the fluctuations and could only reduce the order of accuracy of the method.

\begin{figure}[t!]
	\centering
	\begin{subfigure}[b]{0.495\textwidth}
		\centering
		\includegraphics[width=\textwidth]{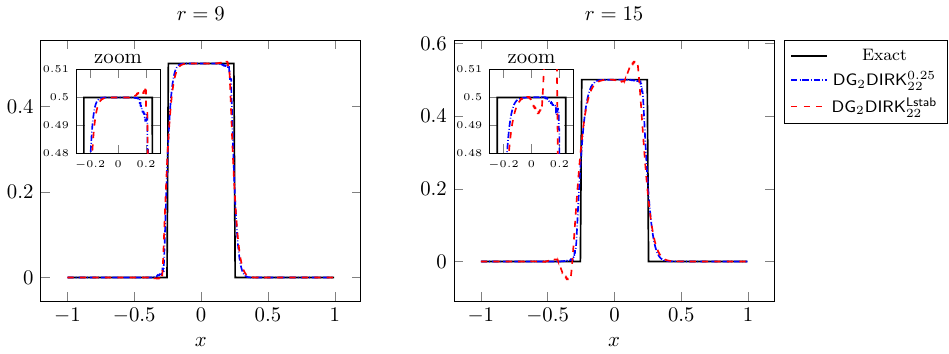}
		\caption{Second-order schemes.\label{fig:li:doubleStep:p1}}
	\end{subfigure}
	\hfill
	\begin{subfigure}[b]{0.495\textwidth}
		\centering
		\includegraphics[width=\textwidth]{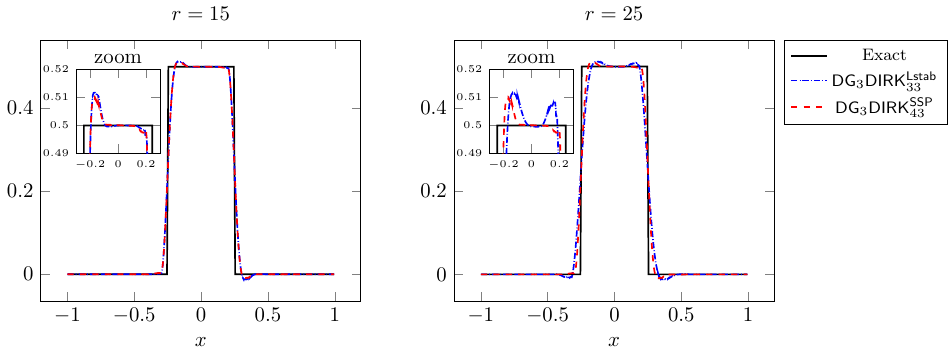}
		\caption{Third-order schemes.\label{fig:li:doubleStep:p2}}
	\end{subfigure}
	\caption{Solutions to the linear advection equation with speed $a=1$ and final time $t=2$. The initial condition is the double-step profile~\eqref{eq:doublestepIC}. In this figure, all the numerical solutions are obtained with $N=400$ cells.\label{fig:li:doubleStep}}
\end{figure}

\begin{figure}[t!]
	\centering
	\begin{subfigure}[b]{0.495\textwidth}
		\centering
		\includegraphics[width=\textwidth]{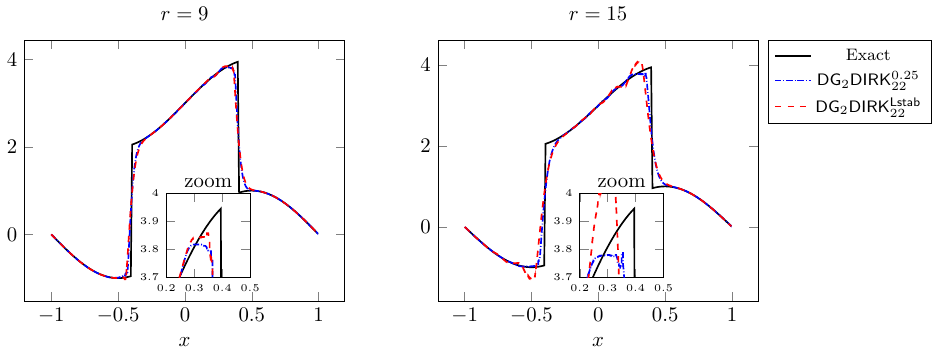}
		\caption{Second-order schemes.\label{fig:li:discSin:p1}}
	\end{subfigure}
	\hfill
	\begin{subfigure}[b]{0.495\textwidth}
		\centering
		\includegraphics[width=\textwidth]{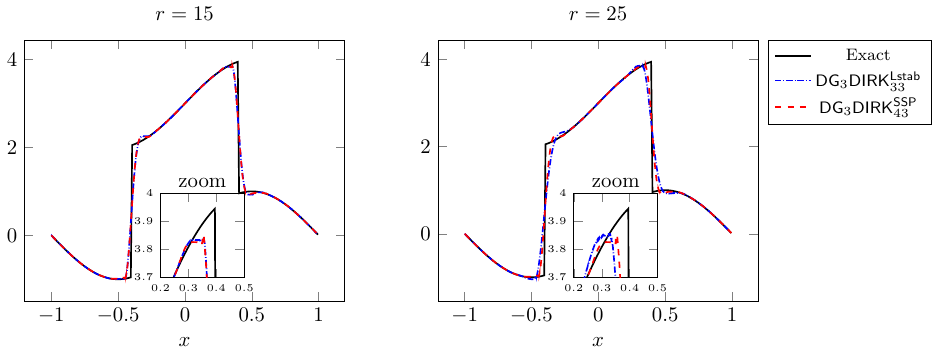}
		\caption{Third-order schemes.\label{fig:li:discSin:p2}}
	\end{subfigure}
	\caption{Solutions to the linear advection equation with speed $a=1$ and final time $t=2$. The initial condition is the discontinuous sinusoidal profile~\eqref{eq:sindiscontIC}. In this figure, all the numerical solutions are obtained with $N=400$ cells.\label{fig:li:discSin}}
\end{figure}

In Figure~\ref{fig:li:doubleStep} and Figure~\ref{fig:li:discSin} we compare the solutions obtained with second- and third-order schemes on the linear advection of the discontinuous initial condition~\eqref{eq:doublestepIC} and~\eqref{eq:sindiscontIC}, respectively. The simulations are run with varying ratios $r$ between the explicit and implicit time-steps. Subfigure~\ref{fig:li:doubleStep:p1} and Subfigure~\ref{fig:li:discSin:p1} show results for second-order schemes, while Subfigure~\ref{fig:li:doubleStep:p2} and Subfigure~\ref{fig:li:discSin:p2} presents third-order schemes. For the second-order schemes, there is a noticeable difference in accuracy between the methods, particularly when $r=15$, where the $\DIRKlstab$ scheme exhibits more pronounced oscillations near the discontinuities due to the larger dispersive error compared to $\DIRKdue$. In contrast, third-order schemes show better agreement with the exact solution, with less oscillation and better resolution of the step profile, in particular when considering the scheme $\DIRKQuattroTre$. This suggests that $\DIRKdue$ and $\DIRKQuattroTre$ are more robust in handling the advection of non-smooth data, especially under higher time-step ratios $r$. The performance difference underscores the importance of choosing an appropriate scheme based on the desired accuracy and the specific characteristics of the problem.

\subsubsection{Burgers' equation}

\begin{figure}[t!]
	\centering
	\begin{subfigure}[b]{0.495\textwidth}
		\centering
		\includegraphics[width=\textwidth]{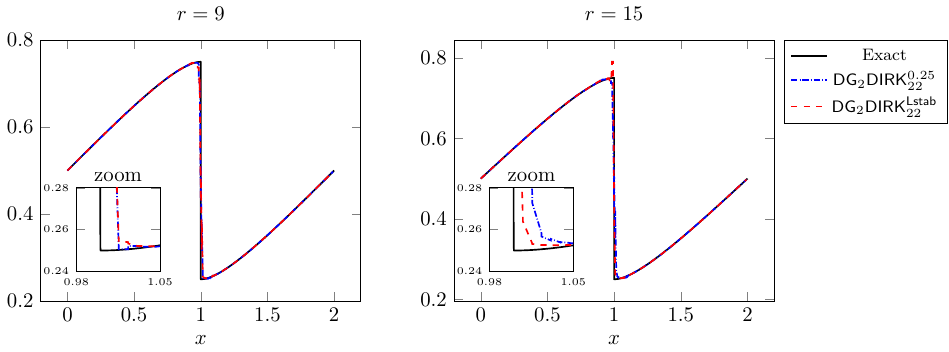}
		\caption{Second-order schemes.\label{fig:burgers:discSin:p1}}
	\end{subfigure}
	\hfill
	\begin{subfigure}[b]{0.495\textwidth}
		\centering
		\includegraphics[width=\textwidth]{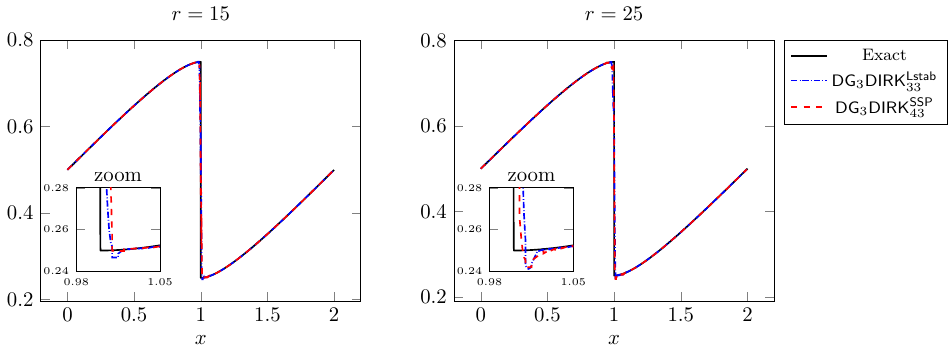}
		\caption{Second-order schemes.\label{fig:burgers:discSin:p2}}
	\end{subfigure}
	\caption{Solutions to the nonlinear Burgers' equation at final time $t=2$ and with smooth initial condition~\eqref{eq:smoothIC}. In this figure, all the numerical solutions are obtained with $N=400$ cells.\label{fig:burgers:discSin:p1:p2}}
\end{figure}

\begin{figure}[t!]
	\centering
	\begin{subfigure}[b]{0.495\textwidth}
		\centering
		\includegraphics[width=\textwidth]{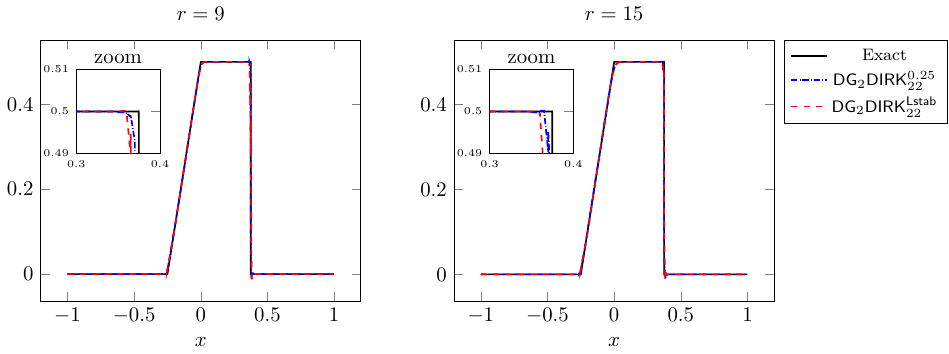}
		\caption{Second-order schemes.\label{fig:burgers:doubleStep:p1}}
	\end{subfigure}
	\hfill
	\begin{subfigure}[b]{0.495\textwidth}
		\centering
		\includegraphics[width=\textwidth]{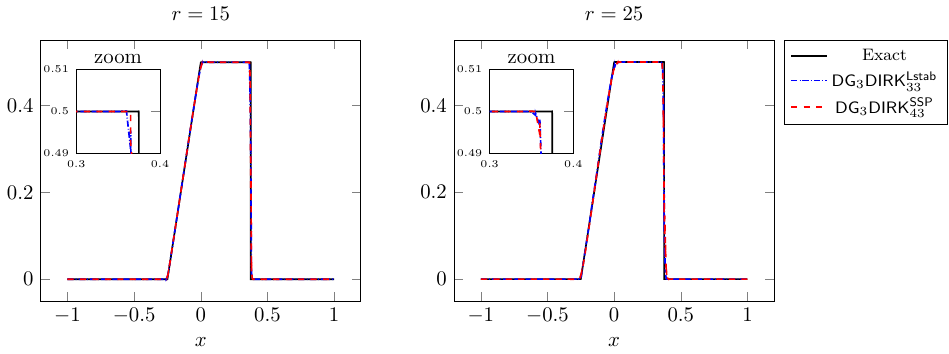}
		\caption{Second-order schemes.\label{fig:burgers:doubleStep:p2}}
	\end{subfigure}
	\caption{Solutions to the nonlinear Burgers' equation at final time $t=0.5$ and with non-smooth initial condition~\eqref{eq:doublestepIC}. In this figure, all the numerical solutions are obtained with $N=400$ cells.\label{fig:burgers:doubleStep:p1:p2}}
\end{figure}

We investigate the behavior of the schemes on the nonlinear Burgers' equation
\begin{equation} \label{eq:burgers}
	u_t + \left( \frac{u}{2} \right)^2_x = 0,
\end{equation}
for different initial conditions. First, we consider the smooth initial condition
\begin{equation} \label{eq:smoothIC}
	u_0(x) = 0.5 - 0.25 \sin(\pi x),
\end{equation}
on the periodic domain $x\in[0,2]$, and up to the final time $t=2$, i.e.~after shock formation. Then, we test the Burgers' equation on the discontinuous initial condition~\eqref{eq:doublestepIC}, on the periodic domain $x\in[-1,1]$, and up to the final time $t=0.5$. Figures~\ref{fig:burgers:discSin:p1:p2} and Figure~\ref{fig:burgers:doubleStep:p1:p2} present the results of different DIRK-DG schemes. \Rtwo{In these numerical tests $\delta$ is chosen as $\delta\approx r/(2p+1)$.} All solutions are computed with a tolerance of $10^{-5}$ for the GMRES solver, and with a maximum number of iterations limited to $100$. On average, the GMRES solver converges within 10 iterations for the second-order schemes and approximately 15 iterations for the third-order schemes, with slightly higher values at initial time steps. A detailed technical study of the iteration behavior is beyond the scope of this work.

\begin{figure}[t!]
	\centering
	\begin{subfigure}[b]{0.49\textwidth}
		\centering
		\includegraphics[width=\textwidth]{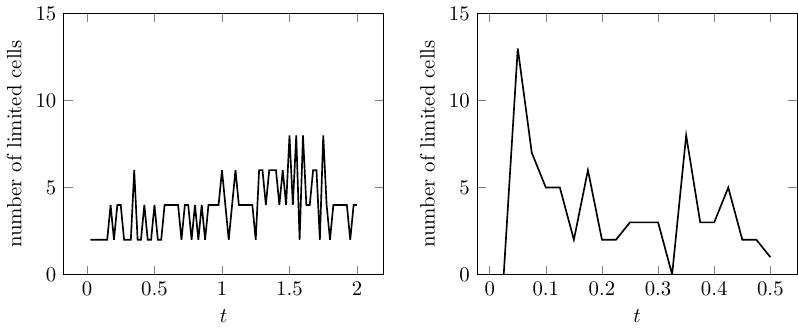}
		\caption{Second-order scheme $\DIRKdue$. Left, case of the initial condition~\eqref{eq:smoothIC}. Right, case of the initial condition~\eqref{eq:doublestepIC}. In both cases, $r=15$.\label{fig:burgers:limiter:p1}}
	\end{subfigure}
	\hfill
	\begin{subfigure}[b]{0.49\textwidth}
		\centering
		\includegraphics[width=\textwidth]{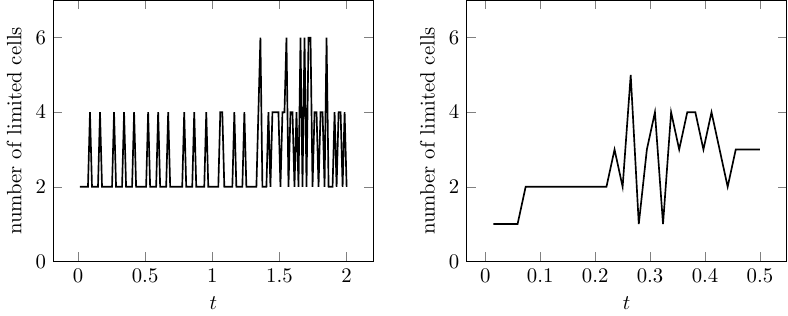}
		\caption{Third-order scheme $\DIRKlstabTre$. Left, case of the initial condition~\eqref{eq:smoothIC}. Right, case of the initial condition~\eqref{eq:doublestepIC}. In both cases, $r=15$\label{fig:burgers:limiter:p2}}
	\end{subfigure}
	\caption{Number of troubling cells detected by the MP limiter~\eqref{eq:limiterMP:pred} computed on the predictor solution. In each time-step, the number of limited cells is chosen as the maximum number of troubling cells detected over the DIRK stages.\label{fig:burgers:limiter}}
\end{figure}

In Figure~\ref{fig:burgers:discSin:p1:p2}, we observe the solutions obtained at the final time $t=2$ using the smooth initial condition described in equation~\eqref{eq:smoothIC}. The numerical solutions are computed using $400$ cells. Subfigure~\ref{fig:burgers:discSin:p1} shows the results for second-order schemes. Subfigure~\ref{fig:burgers:discSin:p2} displays results from the third-order schemes. The scheme $\DIRKdue$ exhibits minimal oscillations and aligns closely with the exact solution, particularly near the peak and smooth regions of the profile. The scheme $\DIRKlstab$ also shows good agreement with the exact solution but introduces slightly more oscillations near the steep gradient.
Both third-order schemes outperform the second-order schemes in terms of accuracy.

Same considerations hold for Figure~\ref{fig:burgers:doubleStep:p1:p2}, where the nonlinear Burgers' equation is solved using a non-smooth initial condition (a double-step profile) at the final time $t=0.5$. Similar to Figure~\ref{fig:burgers:discSin:p1:p2}, the simulations use $400$ cells. The subfigures are organized similarly, with second-order schemes in Subfigure~\ref{fig:burgers:doubleStep:p1} and third-order schemes in Subfigure~\ref{fig:burgers:doubleStep:p2}. Second-order schemes exhibit small oscillations near the shock wave, whereas the third-order schemes show a slightly better resolution, particularly on the rarefaction wave.

In Figure~\ref{fig:burgers:limiter} we illustrate the number of troubling cells detected by the MP limiter during the simulation for the two different initial conditions~\eqref{eq:smoothIC} and~\eqref{eq:doublestepIC}, using the second-order $\DIRKdue$ and third-order $\DIRKlstabTre$ schemes. The MP limiter is applied to the predictor solution, and at each time-step, the number of limited cells represents the maximum number of troubling cells detected across the DIRK stages.
When considering the smooth initial condition, the number of troubling cells starts low and gradually increases as the solution evolves, because the initial condition is smooth, and as time progresses, a shock forms, leading to an increase in the number of cells requiring limiting.

\subsection{Euler system of gas-dynamics}

In this section, we consider several problems based on the one-dimensional nonlinear Euler system for gas-dynamics
\begin{equation} \label{eq:euler:system}
	\partial_t 
	\left( \begin{array}{c}
		\rho \\ \rho v \\ E
	\end{array}\right) +
	\partial_x 
	\left( \begin{array}{c}
		\rho v \\ \rho v^2 + p \\ v(E+p)
	\end{array}\right)  = 0.
\end{equation}
Here, $\rho$, $v$, $p$ and $E$ represent the density, velocity, pressure, and energy per unit volume of an ideal gas, respectively. The equation of state for the gas is given by $E = \frac{p}{\Gamma-1} + \frac12 \rho v^2, $ where $\Gamma = 1.4$. The eigenvalues of the Jacobian of the flux are $\Lambda_1 = v-c$, $\Lambda_2=v$ and $\Lambda_3=v+c$, with the sound speed given by $c=\sqrt{\Gamma p/\rho}$. To focus on the accuracy of material wave, the parameter for numerical viscosity in~\eqref{eq:LFflux} is set to $\alpha(\mathbf{u}^+,\mathbf{u}^-) = \max\{|v^{-}|,|v^{+}|\}$, where $v^{-}$ and $v^{+}$ are the DG reconstructions of the gas velocity at the interfaces.

In the following we consider only the second-order scheme $\DIRKlstab$ and the third-order scheme $\DIRKlstabTre$, which are both L-stable. They will be compared with explicit RKDG schemes and the implicit Quinpi finite-volume scheme of~\cite{Puppo2023,Puppo2024}.

\subsubsection{Convergence test}

We check the order of accuracy of DIRK-DG schemes on the nonlinear Euler system~\eqref{eq:euler:system} by simulating the linear transport of a density perturbation. The initial condition is
\begin{equation} \label{eq:conv:test:ic}
	\left(\rho,v,p\right) = \left(1+0.5\sin(2\pi x), 1, 10^\kappa\right), \quad \kappa\in\mathbb{N}_0,
\end{equation}
and the exact solution is characterized by a traveling wave that at a given time $t$ can be written as
\begin{equation} \label{eq:conv:test:exact}
	\left(\rho,v,p\right) = \left(1+0.5\sin\left(2\pi (x-t)\right), 1, 10^\kappa\right), \quad \kappa\in\mathbb{N}_0.
\end{equation}
Here, $\kappa$ allows as to describe different stiffness regimes. In fact, for this problem, the spectral radius is $\Lambda_{\max} := \max_{i=1,2,3} |\Lambda_i| = 1 + 10^{\kappa/2}\sqrt{2\Gamma}$. The eigenvalues $\Lambda_i$ of the system are preserved during the time evolution. For an explicit RKDG scheme of order $p+1$ the CFL stability condition would require to choose a grid ratio
$$
	\frac{\Delta t_{\text{CFL}}}{ h} \leq \frac{(2p+1)^{-1}}{\Lambda_{\max}} \sim 10^{-\kappa/2} (2p+1)^{-1}.
$$ 
In order to test the accuracy of the implicit scheme, we consider $\kappa=0$ and $\kappa=2$ and run the simulation with the following ratio between the explicit and the implicit time-step
$$
	r = 4 \Lambda_{\max} \approx \begin{cases} 10.69, & \kappa=0 \\ 70.93, & \kappa=2. \end{cases}
$$
\Rtwo{
These choices imply the following implicit time-steps
$$
	\Delta t_{\text{imp}} = \begin{cases}
		3.563, & \kappa=0,\ p=1 \\
		2.138, & \kappa=0,\ p=2 \\
		2.364, & \kappa=2,\ p=1 \\
		1.419, & \kappa=2,\ p=2.
	\end{cases}
$$
}
We run the problem up to the final time $t=0.25$ and use periodic boundary conditions on the domain $[0,1]$.

\begin{table}[th!]
	\caption{Experimental orders of convergence computed on the Euler system~\eqref{eq:euler:system} with initial condition~\eqref{eq:conv:test:ic}. \Rtwo{The errors are computed on the density.}\label{tab:rates:sys}}
	\begin{subtable}[h]{0.45\textwidth}
	\caption{Second-order scheme $\DIRKlstab$.\label{tab:rates:sys:p1}}
	\centering
	\vspace{0.25cm}
	\pgfplotstabletypeset[
	font=\small,
	col sep=comma,
	sci zerofill,
	empty cells with={--},
	every head row/.style={before row={\toprule
			&\multicolumn{2}{c}{$\kappa=0$}
			&\multicolumn{2}{c}{$\kappa=2$}
			\\
		},
		after row=\midrule
	},
	every last row/.style={after row=\bottomrule},
	create on use/rate1/.style={create col/dyadic refinement rate={1}},
	create on use/rate2/.style={create col/dyadic refinement rate={2}},
	columns/0/.style={column name={$N$}},
	columns/1/.style={column name={$L^1$ error},sci e},
	columns/rate1/.style={fixed zerofill,column name={rate}},
	columns/2/.style={column name={$L^1$ error},sci e},
	columns/rate2/.style={fixed zerofill,column name={rate}},
	columns={0,1,rate1,2,rate2},
	]
	{errsys_p1.err}
	\end{subtable}
	\hfill
	\begin{subtable}[h]{0.45\textwidth}
		\caption{Third-order scheme $\DIRKlstabTre$.\label{tab:rates:sys:p2}}
		\centering
		\vspace{0.25cm}
		\pgfplotstabletypeset[
		font=\small,
		col sep=comma,
		sci zerofill,
		empty cells with={--},
		every head row/.style={before row={\toprule
				&\multicolumn{2}{c}{$\kappa=0$}
				&\multicolumn{2}{c}{$\kappa=2$}
				\\
			},
			after row=\midrule
		},
		every last row/.style={after row=\bottomrule},
		create on use/rate1/.style={create col/dyadic refinement rate={1}},
		create on use/rate2/.style={create col/dyadic refinement rate={2}},
		columns/0/.style={column name={$N$}},
		columns/1/.style={column name={$L^1$ error},sci e},
		columns/rate1/.style={fixed zerofill,column name={rate}},
		columns/2/.style={column name={$L^1$ error},sci e},
		columns/rate2/.style={fixed zerofill,column name={rate}},
		columns={0,1,rate1,2,rate2},
		]
		{errsys_p2.err}
	\end{subtable}
\end{table}

Table~\ref{tab:rates:sys} presents the experimental orders of convergence computed for the Euler system under two different schemes: the second-order scheme $\DIRKlstab$ in~\ref{tab:rates:sys:p1} and the third-order scheme $\DIRKlstabTre$ in Table~\ref{tab:rates:sys:p2}. Both tables provide data for different levels of grid refinement and for two different values of $\kappa$, which defines the level of stiffness. \Rtwo{The errors are computed on the density profile.}

We observe that for the second-order scheme, the experimental order of convergence is not exactly as the theoretical value, likely due to the space limiter employed, i.e.~the moment limiter, which uses the minmod function. This function can lead to the clipping of extrema, slightly reducing the expected convergence rate. Also, for $\kappa=2$ the problem is particularly stiff and finer grids can be used to observe the expected order. The third-order scheme, instead, is better able to mitigate the effects of the space limiter and of the stiffness.

\subsubsection{Stiff Riemann problems}

In the following section, we evaluate the performance of the DIRK-DG schemes on stiff Riemann problems for the Euler system~\eqref{eq:euler:system}. Our focus is on scenarios where there is a significant disparity between the propagation speeds of acoustic and material waves. Specifically, we are interested in cases where the ratio $\frac{|v|}{|v|+c} \ll 1, \ \forall\,(x,t)$. This condition highlights the challenge of accurately capturing the interaction between fast-moving acoustic waves and slower material waves.

The acoustic wave speeds impose a strict constraint on the time-step for stability:
$$
\frac{\Delta t_{\text{CFL}}}{ h} \leq \frac{(2p+1)^{-1}}{\max_{x} |v|+c},
$$
where $\Delta t_{\text{CFL}}$ denotes the stability-limited time-step.
However, with implicit schemes, the time-step is no longer constrained by stability requirements. Instead, we can choose the time-step to accurately resolve the slow-moving material wave:
$$
\frac{\Delta t_{\text{acc}}}{ h} = \frac{1}{|v_{\text{cw}}|},
$$
where $v_{\text{cw}}$ represents an estimate of the material wave velocity. The ratio between the stability-limited time-step $\Delta t_{\text{CFL}}$ and the accuracy-driven time-step $\Delta t_{\text{acc}}$ for the material wave is a measure of the stiffness of the problem. This ratio was previously defined by $r$, see~\eqref{eq:def:Courant}, namely
$$
r = \frac{\Delta t_{\text{acc}}}{\Delta t_{\text{CFL}}}.
$$
This number indicates how much the time-step is relaxed by using an implicit scheme, compared to the explicit stability requirement.

We consider the following Riemann problems that have been previously studied in~\cite{Puppo2024}. For reference, Figure 2 in that paper illustrates the exact density solution and the corresponding wave speeds at time $t=1$.

\begin{description}
	\item[Test (a): Non-symmetric expansion problem.]
	\begin{equation} \label{eq:tworar:ic}
		\left(\rho, u, p\right) = \begin{cases}
			(1, -0.15, 1), & x < 0, \\
			(0.5, 0.15, 1), & x \geq 0.
		\end{cases}
	\end{equation}
	
	The exact density solution at time $t = 1$ is characterized by a 1-rarefaction wave propagating with a negative velocity, a 3-rarefaction wave moving with a positive velocity, and a 2-contact wave traveling with a velocity $v_{\text{cw}} = -2.57 \times 10^{-2}$.
	
	For this problem, one has $\max_{(x,t)} |v| + c = 1.82$.
	As a result, an explicit scheme would require a time-step $\Delta t_{\text{CFL}}$ satisfying the stability condition
	$$
	\frac{\Delta t_{\text{CFL}}}{h} \leq \frac{(2p + 1)^{-1}}{1.82} = \frac{0.549}{2p + 1}.
	$$
	
	However, with an implicit scheme, this stability constraint can be bypassed. Specifically, to ensure accuracy in capturing the contact wave, the time-step $\Delta t_{\text{acc}}$ should be chosen such that
	\begin{equation} \label{eq:dt:tworar}
	\frac{\Delta t_{\text{acc}}}{h} = \frac{1}{-2.57 \times 10^{-2}} = 6.66,
	\end{equation}
	which leads to a ratio $r$ between the accuracy-driven time-step and the stability-driven time-step of
	$$
	r = (2p+1)\frac{6.66}{0.549} \approx 12.1(2p + 1).
	$$

	\begin{figure}[t!]
		\centering
		\begin{subfigure}[b]{\textwidth}
			\centering
			\includegraphics[width=\textwidth]{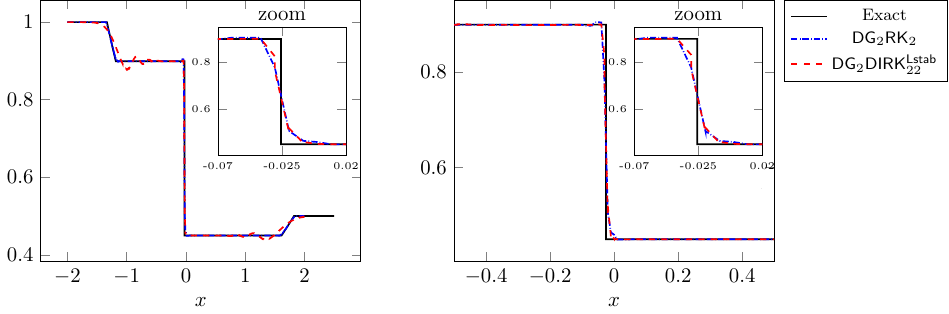}
			\caption{Second-order schemes. Left: density solution with $N=400$ cells. Right: density solution computed on a shorter tube section, around the contact wave, with $N=100$ cells.\label{fig:sys:tworar:p1}}
		\end{subfigure}
		\begin{subfigure}[b]{\textwidth}
			\centering
			\includegraphics[width=\textwidth]{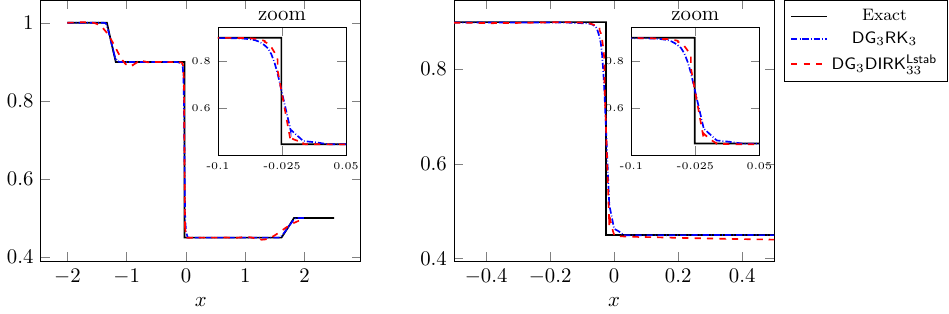}
			\caption{Third-order schemes. Left: density solution with $N=266$ cells. Right: density solution computed on a shorter tube section, around the contact wave, with $N=66$ cells.\label{fig:sys:tworar:p2}}
		\end{subfigure}
		\caption{Solutions to the Euler system of gas-dynamics, with initial condition~\eqref{eq:tworar:ic} at final time $t=1$. Comparison between implicit and explicit schemes.}\label{fig:sys:tworar}
	\end{figure}
	
	\item[Test (b): colliding flow problem.]
	\begin{equation} \label{eq:twoshock:ic}
		\left(\rho,v,p\right) = \begin{cases}
			(1.5,  0.5, 10), & x < 0\\
			(0.5, -0.5, 10), & x \geq 0.
		\end{cases}
	\end{equation}

	The exact density at time $t=1$ is characterized by a 1-shock wave moving with negative velocity, a 3-shock wave moving with positive velocity, and a 2-contact wave with velocity $v_{\text{cw}}=0.13$.
	
	For this problem, one has $\max_{(x,t)} |v|+c = 5.79$. As a result, an explicit scheme would require a time-step $\Delta t_{\text{CFL}}$ satisfying the stability condition
	$$
	\frac{\Delta t_{\text{CFL}}}{ h} \leq \frac{(2p+1)^{-1}}{5.79} = \frac{0.17}{2p+1}.
	$$
	
	However, with an implicit scheme, this stability constraint can be bypassed. Specifically, to ensure accuracy in capturing the contact wave, the time-step $\Delta t_{\text{acc}}$ should be chosen such that
	\begin{equation} \label{eq:dt:twoshock}
	\frac{\Delta t_{\text{acc}}}{ h} = \frac{1}{0.13} = 7.7,
	\end{equation}
	which leads to a ratio $r$ between the accuracy-driven time-step and the stability-driven time-step of
	$$
	r = (2p+1) \frac{7.7}{0.17} \approx 45.3 (2p+1).
	$$
	
	\begin{figure}[t!]
		\centering
		\begin{subfigure}[b]{\textwidth}
			\centering
			\includegraphics[width=\textwidth]{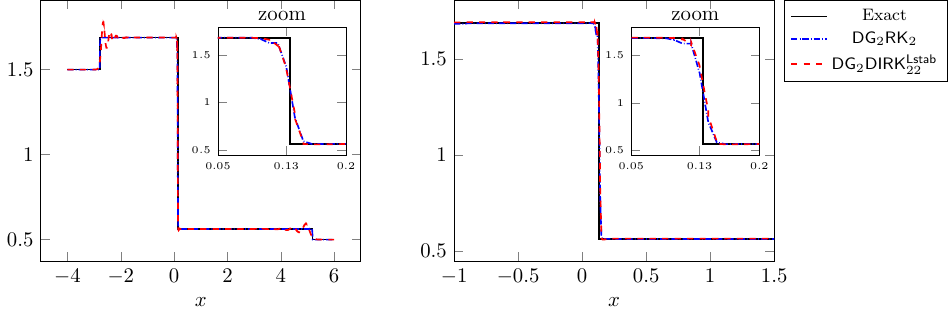}
			\caption{Second-order schemes. Left: density solution with $N=1000$ cells. Right: density solution computed on a shorter tube section, around the contact wave, with $N=250$ cells.\label{fig:sys:twoshock:p1}}
		\end{subfigure}
		\begin{subfigure}[b]{\textwidth}
			\centering
			\includegraphics[width=\textwidth]{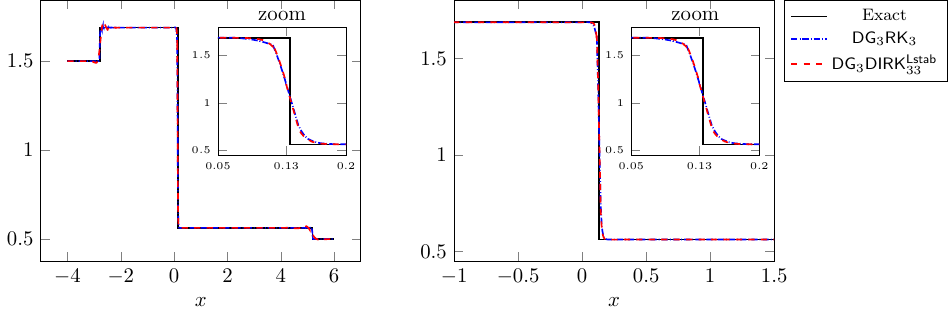}
			\caption{Third-order schemes. Left: density solution with $N=666$ cells. Right: density solution computed on a shorter tube section, around the contact wave, with $N=166$ cells.\label{fig:sys:twoshock:p2}}
		\end{subfigure}
		\caption{Solutions to the Euler system of gas-dynamics, with initial condition~\eqref{eq:twoshock:ic} at final time $t=1$. Comparison between implicit and explicit schemes.}\label{fig:sys:twoshock}
	\end{figure}

%
%
%

\end{description}

Figure~\ref{fig:sys:tworar} and Figure~\ref{fig:sys:twoshock} present a comparison of density solutions to the Euler system of gas dynamics~\eqref{eq:euler:system}, with initial conditions~\eqref{eq:tworar:ic} and~\eqref{eq:twoshock:ic}, respectively, at final time $t=1$. The comparison includes both implicit and explicit schemes for second-order and third-order methods. The explicit schemes use the Heun method for second-order and the SSP Runge-Kutta method for third-order time-stepping. We refer to them as $\DGRKdue$ and $\DGRKtre$, respectively. Explicit schemes are run with $\Delta t_{\text{CFL}}$, whereas implicit schemes employ the time-step $\Delta t_{\text{acc}}$, focused on the material speed.

All solutions are computed using an extension parameter $\delta=1$. A tolerance of $10^{-5}$ was set for the GMRES solver, with a maximum number of iterations limited to $500$. On average, the GMRES solver converges within 45 iterations both for the second-order and the third-order scheme.

The number of cells is chosen such that the degrees of freedom are the same between the second-order and third-order schemes. This ensures a fair comparison in terms of computational accuracy and performance.

The solution is computed on two different computational domains, were the short domain is chosen to be 1/4 of the long domain. The mesh width $h$ is the same on both, long and short, domains.
In the context of stiff problems, where the primary interest lies in accurately reproducing material waves rather than faster-moving acoustic waves, the ability to compute the solution only where the material wave actually is, is important, because it speeds up the computation. In DG schemes, this can be done, because boundary conditions can be easily enforced, due to the compactness of the stencil. In fact, in the short computational domain, one can see that acoustic waves have already left the domain, without interfering with the material wave.

In Figure~\ref{fig:sys:tworar:p1}, the left plot shows the density of the non-symmetric expansion obtained with second-order schemes on $N=400$ cells on the long tube, while the right plot shows the solution on the shorter tube with the same $h$. The left plot of Figure~\ref{fig:sys:twoshock:p1} refers to the density of the colliding flows problem obtained with second-order schemes on $N=1000$ cells on the long tube, while the right plot shows the solution on the shorter tube with the same $h$.

In Figure~\ref{fig:sys:tworar:p2}, we show the solution obtained with the third-order scheme with $N=266$ cells on the long pipe and $N=66$ on the short one. Both grids use the same $h$.

Both explicit and implicit schemes closely follow the exact solutions. The implicit ones exhibit deviations observable in regions near the acoustic waves. Instead, they are comparable to explicit schemes on the contact wave even with a time-step which is much larger. In addition, we observe that, in the case of the short tube section, the acoustic waves can exit the domain without introducing reflections or other boundary-related artifacts that could contaminate the solution.

\begin{figure}[t!]
	\centering
	\begin{subfigure}[b]{\textwidth}
		\centering
		\includegraphics[width=0.45\textwidth]{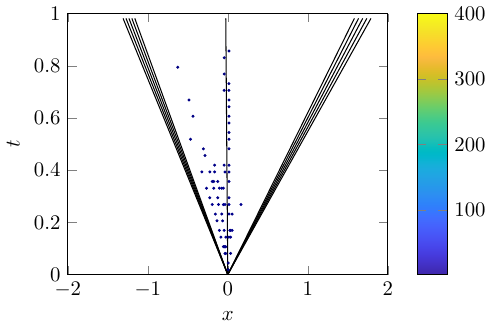}
		\includegraphics[width=0.45\textwidth]{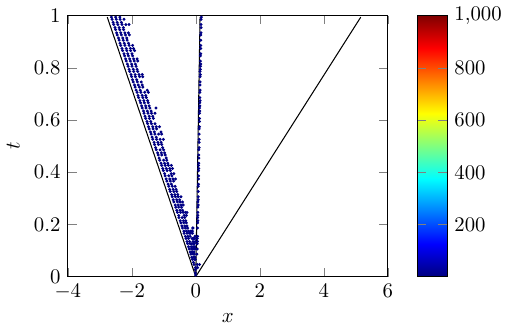}
		\caption{Second-order schemes. Left: density solution \Rtwo{in the $x$-$t$ diagram} of the non-symmetric expansion problem. Right: density solution \Rtwo{in the $x$-$t$ diagram} of the colliding flows problem. \Rtwo{The black lines represent the wave structures.}\label{fig:sys:limiter:p1}}
	\end{subfigure}
	\begin{subfigure}[b]{\textwidth}
		\centering
		\includegraphics[width=0.45\textwidth]{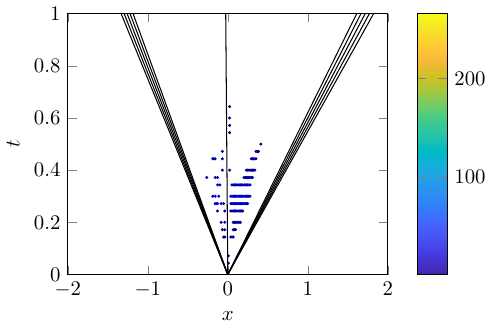}
		\includegraphics[width=0.45\textwidth]{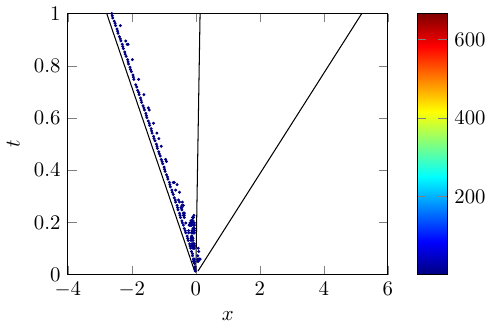}
		\caption{Third-order schemes. Left: density solution \Rtwo{in the $x$-$t$ diagram} of the non-symmetric expansion problem. Right: density solution \Rtwo{in the $x$-$t$ diagram} of the colliding flows problem. \Rtwo{The black lines represent the wave structures.}\label{fig:sys:limiter:p2}}
	\end{subfigure}
	\caption{Structure of the solutions in the $x$-$t$ diagram. Troubling cells are marked at the time levels and locations where detected by the MP limiter computed on the predictor. The colors denote the total number of cells detected at a given time. The limits of the colorbar are set between 1 and the number of cells, with each cell's value corresponding directly to a specific color on the scale.}\label{fig:sys:limiter}
\end{figure}

Figure~\ref{fig:sys:limiter} displays the structure of solutions in the $x$-$t$ diagram for both second-order and third-order schemes, focusing on the Riemann problems we have considered: the non-symmetric expansion problem (left plots) and the colliding flows problem (right plots). The colorbar in each plot represents the total number of troubling cells detected at various time levels, with the color scale ranging from blue (lower number of detected cells) to red (higher number of detected cells). The limits of the colorbar are set between 1 and the number of cells used to discretize the computational domain.

We observe that more troubling cells are detected for the colliding flows problem. 
However, in both case, the color intensity remains mostly within the blue range, implying handling of the problem with few detected troubling cells overall. 
The ratio between the troubling cells and the total number of cells indicates the effectiveness of the numerical scheme. The third-order scheme consistently shows a lower ratio, suggesting a more accurate and stable solution compared to the second-order scheme, particularly on the colliding flows and in regions of complex interactions. This is consistent with the diffusion-dispersion study, from which the high dispersion of second-order schemes was apparent.

\subsubsection{Comparison with the third-order Quinpi finite-volume scheme}

\begin{figure}[t!]
	\centering
	\begin{subfigure}[b]{\textwidth}
		\centering
		\includegraphics[width=\textwidth]{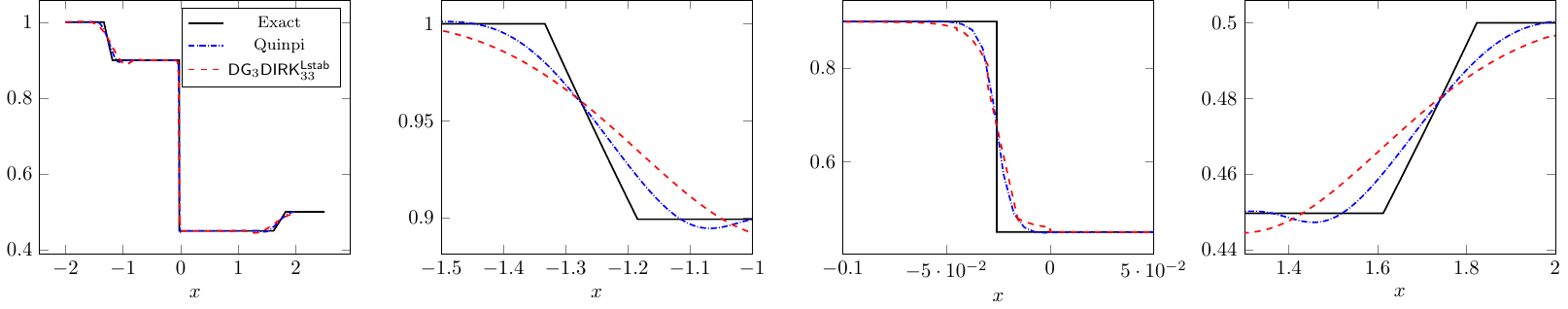}
		\caption{Non-symmetric expansion problem. The left-most panel shows the density solution. The other panels show zooms on the $1$-rarefaction, the contact and $2$-rarefaction wave, respectively.\label{fig:sys:tworar:quinpi}}
	\end{subfigure}
	\begin{subfigure}[b]{\textwidth}
		\centering
		\includegraphics[width=\textwidth]{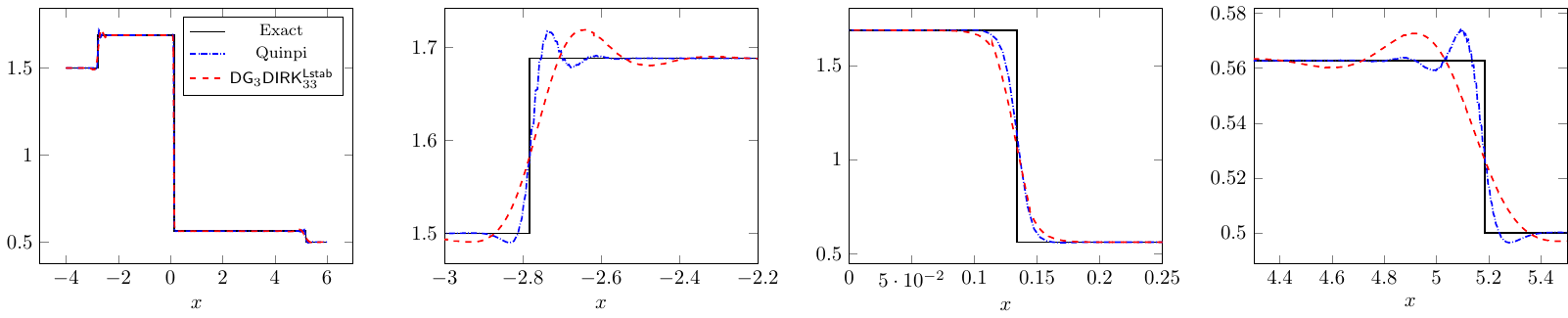}
		\caption{Colliding flow problem. The left-most panel shows the density solution. The other panels show zooms on the $1$-shock the contact and $2$-shock wave, respectively\label{fig:sys:twoshock:quinpi}}
	\end{subfigure}
	\caption{Comparison between the third-order DIRK-DG implicit scheme and the third-order finite volume implicit scheme from~\cite{Puppo2024}. For the non-symmetric expansion problem, N=266 cells are used for the $\DIRKlstabTre$ scheme and $N=800$ cells for Quinpi. For the colliding flow problem, $N=666$ cells are used for the $\DIRKlstabTre$ scheme and $N=2000$ cells for Quinpi. The number of cells is chosen to ensure that both schemes are computed with \Rtwo{comparable number of degrees of freedom}.}\label{fig:sys:quinpi}
\end{figure}

Figure~\ref{fig:sys:quinpi} provides a detailed comparison between the third-order DIRK-DG implicit scheme and the third-order finite volume implicit scheme from~\cite{Puppo2024} on the non-symmetric expansion problem (Figure~\ref{fig:sys:tworar:quinpi}) and the colliding flow problem (Figure~\ref{fig:sys:twoshock:quinpi}). For the non-symmetric expansion, the simulations were carried out using $N=266$ cells for the $\DIRKlstabTre$ scheme and $N=800$ cells for Quinpi. In the colliding flow scenario, $N=666$ cells were used for the $\DIRKlstabTre$ scheme and $N=2000$ cells for Quinpi. The choice of cell numbers was made to ensure that both schemes operated with \Rtwo{comparable number of degrees of freedom}, facilitating a fair comparison between them.

Both schemes were run with a fixed ratio between time-step and cell size, given by~\eqref{eq:dt:tworar} for the non-symmetric expansion and by~\eqref{eq:dt:twoshock} for the colliding flows. However, a significant distinction lies in the fact that the stability condition for explicit DG schemes is constrained by the degree of the polynomial used in the discretization, unlike finite-volume schemes where the stability condition is generally less restrictive. This constraint means that explicit time-stepping methods for DG schemes require smaller time-steps to maintain stability, making the ratio between explicit and implicit time-steps much larger for DG schemes when compared to finite-volume schemes. As a result, the DIRK-DG scheme can take much larger time-steps in its implicit form, enhancing its efficiency and stability in challenging scenarios.

The zoomed-in views in Figure~\ref{fig:sys:quinpi} provide a closer look at the small-scale differences between the methods. The figure shows that both the DIRK-DG and Quinpi schemes successfully capture the sharp discontinuity of the exact solution across the material wave of both problems. Instead, they show deviations across the sound waves as expected. On all the waves, the Quinpi scheme looks less diffusive than the DIRK-DG scheme. This behavior will be further investigated, but clearly DG has an edge in the treatment of boundary conditions.

\section{Conclusions} \label{sec:conclusion}

This work aimed at investigating implicit discontinuous Galerkin schemes characterized by low dissipation and dispersion. Through comprehensive Fourier analysis, we identified the optimal combination of DIRK and DG schemes, formulating dispersion and dissipation errors analytically for second-order schemes and numerically for third-order schemes. A significant contribution of this work is the proposal of a method, inspired by the Quinpi approach from~\cite{Puppo2024}, to mitigate the challenges posed by the high nonlinearity introduced by space limiters in implicit methods.

Our results highlight the strengths of the newly developed implicit discontinuous Galerkin schemes. Specifically, these methods demonstrate the capability to accurately capture slow material waves in stiff problems, effectively overcoming the stringent stability constraints, based on fast acoustic waves, that typically limit explicit schemes. This advantage is particularly evident in scenarios where the phenomena of interest propagate at speeds much slower than the fastest wave speeds, allowing the implicit schemes to operate with time-steps much larger than those required by explicit methods, while still achieving better accuracy on the phenomenon of interest. \Rone{At the same time, our results confirm that this advantage comes at the cost of reduced accuracy for fast acoustic waves when large time steps are employed, a limitation that is intrinsic to implicit time integration. Such a trade-off is acceptable in applications where the phenomena of interest are dominated by slow dynamics.}

However, the study also identifies areas where improvements are needed. The solution of the nonlinear systems, a critical component of the implicit discontinuous Galerkin schemes, could benefit from more efficient strategies. Implementing preconditioning techniques to accelerate the GMRES solver is one potential avenue for enhancing computational efficiency. Nevertheless, this preconditioning has to be build on the knowledge of the low-order solution only.

Importantly, this implicit approach has not been tailored to a specific problem, suggesting that it could have broad applicability across a range of stiff problems, particularly those involving slow-moving material waves. The ability of implicit schemes to handle these challenging scenarios with greater accuracy and larger time-steps makes them a promising tool for future applications.

Despite the promising findings, there are some limitations to our study. The Fourier analysis for third-order schemes, in particular, could not be fully addressed analytically, highlighting a potential area for further investigation.

Looking ahead, future work should focus on extending these methods to two-dimensional problems, where the complexity of the wave interactions increases. Additionally, developing more efficient techniques for solving the nonlinear systems inherent in implicit schemes will be crucial to making these methods more practical and widely applicable.

In conclusion, the implicit discontinuous Galerkin schemes developed in this study represent a significant step forward in the numerical treatment of stiff problems, particularly those involving slow material waves. With further refinement and extension, these methods have the potential to make a substantial impact in various scientific and engineering applications.

\section*{Acknowledgments and funding}

This work was partly supported by MUR (Ministry of University and Research) under the PRIN-2022 project on ``High order structure-preserving semi-implicit schemes
for hyperbolic equations'' (number 2022JH87B4).\\
This work was carried out within the Ateneo Sapienza Project 2023 ``Modeling, numerical treatment of hyperbolic equations and optimal control problems'' funded by Sapienza University of Rome.\\
GV acknowledge support from PNRR-MUR (Ministry of University and Research) project ``Italian Research Center on High Performance Computing, Big Data and Quantum Computing''.\\
The authors are members of the INdAM Research National Group of Scientific Computing (INdAM-GNCS).

\section*{Conflict of interest}
There is no conflict of interest in the present work.

\section*{Data Availability Statement}
Data sharing not applicable to this article as no datasets were generated or analysed during the current study.

\end{document}